%% file: main.tex
\documentclass[a4paper,fleqn]{cas-dc}

\usepackage[numbers,sort&compress]{natbib}
\usepackage{amsmath,amssymb}

\newtheorem{theorem}{Theorem}[section]
\newtheorem{lemma}[theorem]{Lemma}
\newtheorem{proposition}[theorem]{Proposition}
\newtheorem{corollary}[theorem]{Corollary}
\newtheorem{definition}[theorem]{Definition}
\newtheorem{example}[theorem]{Example}
\newtheorem{remark}[theorem]{Remark}
\newtheorem{question}[theorem]{Question}
\newenvironment{proof}{\par\noindent\textbf{Proof.}}{\hfill$\square$\par}

\newcommand{\OO}{\mathcal{O}}
\newcommand{\KK}{\mathcal{K}}
\newcommand{\LL}{\mathcal{L}}
\newcommand{\NN}{\mathbb{N}}
\newcommand{\QQ}{\mathbb{Q}}
\newcommand{\RR}{\mathbb{R}}
\newcommand{\Sier}{\mathbb{S}}
\newcommand{\W}{\mathsf{W}}
\newcommand{\Ww}{\mathsf{W}_{w}}
\newcommand{\Wi}{\mathsf{W}_{i}}
\newcommand{\Win}{\mathsf{W}_{in}}
\newcommand{\cl}{\operatorname{cl}}
\newcommand{\sat}{\operatorname{sat}}
\newcommand{\id}{\mathrm{id}}
\newcommand{\ph}{\mathrm{P}_{\!H}}
\newcommand{\ps}{\mathrm{P}_{\!S}}
\newcommand{\Ck}{C_{k}}
\newcommand{\Ok}{\OO_{k}}
\newcommand{\OL}{\OO_{L}}

\begin{document}
\let\WriteBookmarks\relax
\def\floatpagepagefraction{1}
\def\textpagefraction{.001}

\shorttitle{On Wilker-type properties}
\shortauthors{W. Wang and J. J. Tang}

\title[mode=title]{On Wilker-Type Properties}

\tnotemark[1]
\tnotetext[1]{This research is supported by by the Scientific Research Program of the Tianjin Municipal Education
Commission (No.~2023KJ281).}

\author[1]{Wu Wang}
\ead{wangwu@alu.scu.edu}

\author[1]{Jiajun Tang}
\ead{tt20982025@163.com}
\cormark[1]
\affiliation[1]{organization={Zhonghuan Information College, Tianjin University of Technology},
            city={Tianjin},
            postcode={300380},
            country={China}}

\cortext[1]{Corresponding author}

\begin{abstract}
Wilker-type properties record how a compact set can be split along a binary open cover, and they are precisely the continuity conditions for the union operation on open lattices, governing the topology of power and function spaces. We give a unified treatment of the four properties $\W$, $\Ww$, $\Wi$ and $\Win$. We prove that $\Wi$ and $\Win$ are hereditary for arbitrary subspaces, whereas $\W$ and $\Ww$ are hereditary only for closed subspaces; all four are preserved by continuous retracts, and $\Ww$, $\Wi$, $\Win$ are finitely additive over topological sums. For KC-spaces $\Ww$ and $\Wi$ coincide, and spaces whose compact subspaces are Hausdorff belong to $\Win$. Using Kov\'ar's compact KC-space $M_0$ and its variant $\QQ^{+}$, we obtain countably infinite compact KC-spaces in $\Wi\setminus\W$, showing $\Wi$ is strictly weaker than $\W$. Conversely, we construct compact $T_1$ spaces whose product fails $\Ww$, proving the classes $\Ww$ and $\Wi$ are not reflective. We further show Vietoris constructions need not preserve $\Wi$, and prove that $\Ok(M_0)\cong \Ck(M_0,\Sier)$ is not sober and fails property~$S$. Open problems are collected at the end.
\end{abstract}


\maketitle

\input{sec1}
\input{sec2}
\input{sec3}
\input{sec4}
\input{sec5}

\bibliographystyle{elsarticle-num}
\bibliography{wilker}

\end{document}

%% file: sec1.tex
\section{Introduction}\label{sec:intro}

Power spaces form one of the central objects of domain theory, because they supply semantic models
for programming languages endowed with nondeterministic choice
\citep{smyth1983power,gierz2003continuous,chen2024upper}. Classically, the upper power space of a
space $X$ is constructed from the nonempty compact saturated subsets of $X$, and it is again a domain
whenever $X$ is a continuous domain; the construction is then well behaved both algebraically and
order-theoretically. When $X$ is an arbitrary topological space, however, compact saturated sets no
longer retain these properties, and the hypotheses underlying the construction must be identified
individually. The hypotheses that prove to be relevant are refinement
properties of compact sets with respect to open covers, and they can be formulated without any
reference to order theory.

The prototype of such a property is the following question. Let $X$ be a topological space, let
$K\subseteq X$ be compact and let $U_1,U_2\subseteq X$ be open with $K\subseteq U_1\cup U_2$. In
what sense can $K$ be decomposed into two parts $L_1\subseteq U_1$ and $L_2\subseteq U_2$ subject
to a prescribed smallness requirement? If the pieces are required to be compact, one speaks of the
\emph{Wilker property} $\W$; if they are only required to be Lindel\"of, of the \emph{weak Wilker
property} $\Ww$. In both properties the decomposition is effected by subsets of the ambient space.
If, in addition,
the pieces are required to lie inside $K$, one obtains the \emph{inner weak Wilker property} $\Wi$
and the \emph{inner compact refinement property} $\Win$; in that case the requirement
$K\subseteq L_1\cup L_2$ forces $K=L_1\cup L_2$, so that the decomposition is a genuine
decomposition of $K$ itself. The terminology goes back to the study of Wilker spaces and of the
sobriety of function spaces carried out in~\cite{xu2026wilker}, where the four properties were
introduced and compared for the first time.

The interest in these properties is not merely combinatorial. They are equivalent to continuity
statements for the binary union operation on lattices of open sets, namely
\[
  \mu\colon \OL(X)\times \OL(X)\longrightarrow \Ok(X),\qquad \mu(U,V)=U\cup V,
\]
where $\OL(X)$ and $\Ok(X)$ carry the topologies generated by the sets
$\Phi(L)=\{U\in \OO(X):L\subseteq U\}$ for $L$ Lindel\"of and for $L$ compact respectively. The weak
Wilker property is exactly the continuity of $\mu$, and the inner weak Wilker property is exactly
the continuity of all the restrictions $\mu_K$ obtained by replacing $X$ with a compact subspace $K$
\citep{xu2026wilker}. Continuity of such operations is in turn the technical heart of the theory of
power spaces \citep{michael1951topologies,smyth1983power,gierz2003continuous}, and it is also the
key to the sobriety of function spaces with the compact-open topology: the property~$S$ introduced
in~\cite{xu2022sobriety}, which characterizes sobriety of $\Ck(X,Y)$ for sober $Y$ under a mild
order-theoretic hypothesis, is a closure condition that can only be verified through precisely this
kind of refinement argument.

Three groups of questions arising from this circle of ideas remain open, and they constitute the
object of the present paper.

\begin{enumerate}
\item \emph{Permanence.} Which of the four properties are inherited by subspaces, by closed
subspaces, by retracts, by topological sums? A property of compact sets that is not inherited by
closed subspaces is of limited utility, since closed subspaces arise naturally in the
constructions of power spaces.

\item \emph{Separation.} Are the four properties distinct? In particular, is $\Wi$ strictly weaker
than $\W$, that is, is there a space in which every compact set can be decomposed into Lindel\"of
pieces lying inside it, but not into compact pieces? Without such examples the hierarchy collapses
and the inner conditions lose their meaning.

\item \emph{Preservation and reflection.} Do the properties survive the passage to products, to the
power spaces $\ps(X)$ and $\ph(X)$, and to the function space $\Ck(X,\Sier)$? Do they entail
sobriety? Questions of this type determine whether the classes $\W$, $\Ww$, $\Wi$ can serve as a
working category of spaces, for instance whether they are reflective in the category of $T_0$
spaces.
\end{enumerate}

Our answer to the first group is complete for subspaces, retracts and sums. The inner properties
$\Wi$ and $\Win$ are hereditary for arbitrary subspaces, while $\W$ and $\Ww$ are hereditary for
closed subspaces only, and we show by example that the latter cannot be improved: there is a
locally supercompact space in $\W\setminus\Wi$, so $\Ww$ fails to be hereditary for arbitrary
subspaces. All four properties are preserved by continuous retracts, and $\Ww$, $\Wi$, $\Win$ are
finitely additive over topological sums. A useful criterion is that, for KC-spaces, the weak and the
inner weak Wilker properties coincide, and that a space all of whose compact subspaces are Hausdorff
belongs to $\Win$; both criteria reduce the verification of the inner conditions to compact
Hausdorff pieces, where normality is available.

For the second group we use Kov\'ar's countably infinite compact Fr\'echet--Urysohn hyperconnected
KC-space $M_0$ \citep{kovar2005maximal}. Since $M_0$ is countable, every subset of it is Lindel\"of,
which gives $\Wi$; hyperconnectedness together with the KC property forbids a decomposition into two
proper closed pieces, which excludes $\W$. The same argument applies to the space
$\QQ^{+}=\QQ\cup\{\infty\}$ obtained by adding to $\QQ$ a point whose neighbourhoods are the
complements of compact subsets of $\QQ$. Thus $\Wi$ is strictly weaker than $\W$, and the two
properties are separated by countably infinite compact KC-spaces.

The third group presents the principal difficulties, and our results in this direction are
predominantly negative.
Compactness is not preserved by products in this context, and we exploit this in two ways. First,
following the cofinite modification of compact $T_1$ spaces \citep{xu2026cofinite}, we construct
compact $T_1$ spaces $P$ and $Q$, modelled on a grid of cardinality $\aleph_2$ augmented by two
points, such
that every subspace of $P$ and of $Q$ lies in $\Win$, while $P\times Q$ fails even $\Ww$; the
obstruction arises from a cardinality argument showing that no two Lindel\"of pieces cover the grid.
Since $P\sqcup Q$ then belongs to $\Win$ but its square does not, the full subcategories determined
by $\Ww$ and by $\Wi$ are not reflective in the category of $T_0$ spaces. Second, we combine the
same grid with a hereditarily compact $T_1$ topology and with a one-point compactification, and
obtain a hereditarily compact $T_1$ space $E$ and a compact Hausdorff space $Y$ such that
$E\times Y\in\Ww\setminus\Wi$; this shows that the product theorem for $\Ww$ is sharp, and that
$\Wi$ is not preserved by products even when both factors satisfy strong compactness conditions. We also show that
the Vietoris constructions do not preserve $\Wi$: the upper Vietoris space $\ps(Z)$ and the lower
Vietoris space $\ph(Z)$ of a compact $T_1$ space $Z\in\Win$ may both fail $\Ww$, and this failure
persists over locally compact Hausdorff spaces. Finally, we prove that
$\Ok(M_0)\cong \Ck(M_0,\Sier)$ is not sober, and that $\Ck(M_0,\Sier)$ fails property~$S$; this gives
a concrete counterexample in the theory of sobriety of function spaces, and it shows that the
sobriety hypothesis in the characterization given in~\cite{xu2022sobriety} cannot be omitted.

Taken together, the results show that the four Wilker-type properties behave like a genuine
hierarchy of topological invariants: they are stable under the constructions that preserve
compactness, they collapse over KC-spaces, and they are separated by explicit countable compact
examples. The counterexamples indicate that the obstructions are of a cardinality-theoretic rather
than of an order-theoretic nature, and this is the point of view adopted throughout.

%% file: sec2.tex
\section{Preliminaries}\label{sec:prelim}

We collect in this section the notions and the elementary facts that are used repeatedly in the
sequel. For the standard material on topology we refer to Engelking~\cite{engelking1989general}, and
for domain theory and the theory of continuous lattices to Gierz et al.~\cite{gierz2003continuous}
and Abramsky and Jung~\cite{abramsky1994domain}. Throughout, a space is a topological space and maps
are continuous; we write $\OO(X)$ for the family of all open subsets of $X$ and $\cl A$,
$\operatorname{int}A$ for the closure and the interior of $A\subseteq X$.

\subsection{Separation, specialization and saturation}

\begin{definition}\label{def:T0}
A space $X$ is $T_0$ if distinct points have distinct closures, and $T_1$ if every singleton is
closed. On a $T_0$ space the \emph{specialization order} is defined by
\[
  x\le y \quad\Longleftrightarrow\quad x\in \cl\{y\}.
\]
It is a partial order, and every open set is upward closed with respect to it. A subset
$A\subseteq X$ is \emph{saturated} if $A=\uparrow A$, where
$\uparrow A=\{y:\exists a\in A,\ a\le y\}$; equivalently, $A$ is an intersection of open sets. The
\emph{saturation} of an arbitrary subset $A\subseteq X$ is
\[
  \sat A=\bigcap\{U\in\OO(X): A\subseteq U\},
\]
so that $\sat A$ is the smallest saturated set containing $A$.
\end{definition}

Thus a subset of a $T_0$ space is saturated exactly when it is an intersection of open sets, and
$\sat A$ is saturated for every $A$. In a $T_1$ space the specialization order is discrete and all
subsets are saturated, so that saturation carries no information; it becomes meaningful only for
$T_0$ spaces that are not $T_1$, which is precisely the situation in the study of power spaces.

\begin{definition}\label{def:sep}
Let $X$ be a space.
\begin{enumerate}
\item $X$ is \emph{Hausdorff} if distinct points have disjoint open neighbourhoods, and
\emph{normal} if disjoint closed sets have disjoint open neighbourhoods.
\item $X$ is \emph{KC} if every compact subset of $X$ is closed.
\item $X$ is \emph{hyperconnected} if any two nonempty open subsets of $X$ intersect.
\item $X$ is \emph{Fr\'echet--Urysohn} if for every $A\subseteq X$ and every $x\in \cl A$ there is a
sequence $(a_n)$ in $A$ converging to $x$.
\end{enumerate}
\end{definition}

Every Hausdorff space is KC, and every KC space is $T_1$; neither implication reverses. In a
hyperconnected space the complement of a proper closed set is a nonempty open set, so the union of
two proper closed sets is proper. This remark is used in \S\ref{sec:heredity} to exclude the Wilker
property.

\subsection{Compact and Lindel\"of subsets}

\begin{definition}\label{def:compact}
A subset $K$ of a space $X$ is \emph{compact} if every cover of $K$ by open subsets of $X$ has a
finite subcover, and \emph{Lindel\"of} if every such cover has a countable subcover. We write
$\KK(X)$ for the family of all nonempty compact subsets of $X$ and $\LL(X)$ for the family of all
nonempty Lindel\"of subsets of $X$. For $A\subseteq X$ we put
$\LL(A)=\{L\in\LL(X): L\subseteq A\}$ and $\KK(A)=\{K\in\KK(X):K\subseteq A\}$, so that $\LL(A)$ and
$\KK(A)$ are the families of nonempty Lindel\"of, respectively compact, subsets of the subspace $A$.
\end{definition}

Compactness and Lindel\"ofness are intrinsic properties of the subspace $K$, but it is convenient to
test them with covers by open sets of the ambient space, as in Definition~\ref{def:compact}; the two
readings agree. The following permanence facts will be used without further mention.

\begin{lemma}\label{lem:perm}
Let $X$ be a space and let $B\subseteq X$.
\begin{enumerate}
\item The union of finitely many compact subsets of $X$ is compact, and the union of countably many
Lindel\"of subsets of $X$ is Lindel\"of.
\item If $K\subseteq B$ is compact in the subspace $B$, then $K$ is compact in $X$; conversely, if
$K\in\KK(X)$ and $K\subseteq B$, then $K$ is compact in the subspace $B$. The same two statements
hold with ``compact'' replaced by ``Lindel\"of''.
\item If $B$ is closed in $X$ and $L\in\LL(X)$, then $L\cap B$ is Lindel\"of as a subspace of $L$,
hence $L\cap B\in\LL(X)$ whenever $L\cap B\neq\emptyset$.
\item If $f\colon X\to Y$ is continuous, then $f(K)\in\KK(Y)$ for every $K\in\KK(X)$ and
$f(L)\in\LL(Y)$ for every $L\in\LL(X)$.
\item $\sat L$ is Lindel\"of for every $L\in\LL(X)$, and $\Phi(L)=\Phi(\sat L)$ for every
$L\subseteq X$, where $\Phi(E)=\{U\in\OO(X):E\subseteq U\}$.
\end{enumerate}
\end{lemma}

\begin{proof}
(i) is standard. For (ii), let $\mathcal{W}$ be a cover of $K$ by sets open in $B$, say
$W=B\cap U_W$ with $U_W\in\OO(X)$. Since $K$ is compact in $B$, finitely many members cover $K$, and
these are open in $X$; this gives the first claim. Conversely, if $\mathcal{V}$ is a cover of $K$ by
sets open in $B$ and $K$ is compact in $X$, choose for each $V\in\mathcal{V}$ an open set $U_V$ of
$X$ with $V=B\cap U_V$; the family $\{U_V\}$ covers $K$, so finitely many of them, say
$U_{V_1},\dots,U_{V_n}$, cover $K$, and then $V_1,\dots,V_n$ cover $K$. The Lindel\"of version is
proved in the same way, with ``finite'' replaced by ``countable''.

(iii) The set $L\cap B$ is closed in the subspace $L$, because $B$ is closed in $X$; a closed subset
of a Lindel\"of space is Lindel\"of, and Lindel\"ofness is intrinsic by (ii).

(iv) is the standard fact that continuous images of compact, respectively Lindel\"of, spaces are
compact, respectively Lindel\"of.

(v) Let $\mathcal{G}$ be a cover of $\sat L$ by open sets. It covers $L$, so countably many members
$G_1,G_2,\dots$ cover $L$; their union is an open set containing $L$, hence it contains $\sat L$ by
the definition of saturation. Thus $G_1,G_2,\dots$ cover $\sat L$. Finally, $L\subseteq \sat L$
gives $\Phi(\sat L)\subseteq \Phi(L)$, and if $L\subseteq U$ with $U$ open, then $\sat L\subseteq U$,
so $\Phi(L)\subseteq \Phi(\sat L)$.
\end{proof}

\subsection{Power spaces and function spaces}

For a space $X$ we consider the hyperspaces of closed, respectively of compact saturated, subsets of
$X$:
\[
  \begin{gathered}
    \ph(X)=\{F\subseteq X: F=\cl F,\ F\neq\emptyset\},\\
    \ps(X)=\{C\in\KK(X): C=\sat C\}.
  \end{gathered}
\]
For $U\in\OO(X)$ put
\[
  \begin{gathered}
    \langle U\rangle=\{F\in\ph(X): F\cap U\neq\emptyset\},\\
    \Box U=\{C\in\ps(X): C\subseteq U\}.
  \end{gathered}
\]

\begin{definition}\label{def:hyperspace}
The family $\{\langle U\rangle: U\in\OO(X)\}$ is a subbase for the \emph{lower Vietoris topology}
on $\ph(X)$, and the family $\{\Box U: U\in\OO(X)\}$ is a base for the \emph{upper Vietoris
topology} on $\ps(X)$ \citep{michael1951topologies,gierz2003continuous}. The resulting spaces are
denoted by $\ph(X)$ and $\ps(X)$ as well. For spaces $X,Y$ we write $\Ck(X,Y)$ for the set of all
continuous maps $f\colon X\to Y$ endowed with the compact-open topology \citep{isbell1975function}, a subbase of which is given
by
\[
  \begin{gathered}
    [K,V]=\{f\in C(X,Y): f(K)\subseteq V\},\\
    K\in\KK(X),\quad V\in\OO(Y).
  \end{gathered}
\]
The Sierpi\'nski space is $\Sier=\{0,1\}$ with $\{1\}$ as its only nonempty proper open set; it is
the space of truth values, and $\OO(X)$ is in bijection with $C(X,\Sier)$ through the characteristic
functions $\chi_U$.
\end{definition}

Note that the members of $\ps(X)$ are compact by definition, a property on which the utility of the
upper Vietoris topology rests: the sets $\Box U$ form a base encoding the relation of inclusion.

\begin{definition}\label{def:sober}
A nonempty subset $F$ of a space $T$ is \emph{irreducible} if $F\cap G_i\neq\emptyset$ for
$i=1,2$ with $G_1,G_2\in\OO(T)$ implies $F\cap G_1\cap G_2\neq\emptyset$. A \emph{system of
irreducible sets} is a family $H$ that assigns to each space $T$ a collection $H(T)$ of irreducible
subsets of $T$ containing all singletons and closed under continuous images. A point $t\in T$ is a
\emph{generic point} of $F\in H(T)$ if $\cl F=\cl\{t\}$, and $T$ is \emph{$H$-sober} if every
$F\in H(T)$ has a unique generic point. Taking for $H$ the system of all irreducible sets one
recovers the usual notion of sobriety \citep{johnstone1982stone}.
\end{definition}

\begin{definition}\label{def:propS}
Let $X,Y$ be spaces and let $\tau$ be a topology on $C(X,Y)$ whose specialization order coincides
with the pointwise order. Following~\cite{xu2022sobriety}, we say that $\tau$ satisfies
\emph{property $S$} if for every $F\in H(C(X,Y),\tau)$ and every $g\in C(X,Y)$ the identities
\begin{equation}\label{eq:propS}
  g^{-1}(V)=\bigcup_{f\in F} f^{-1}(V)\quad\text{for all } V\in\OO(Y)
\end{equation}
imply $g\in \cl_\tau F$.
\end{definition}

Thus property $S$ is a closure test: a function whose open preimages are unions of the corresponding
preimages of the members of an irreducible family must lie in the closure of that family. For sober
$Y$ and for topologies finer than the pointwise one, this test is equivalent to sobriety of
$(C(X,Y),\tau)$ \citep{xu2022sobriety}; we return to this equivalence in \S\ref{sec:counter}.

\subsection{Core compactness and supercompactness}

\begin{definition}\label{def:corecompact}
Let $U,V$ be open subsets of a space $X$. We say that $V$ is \emph{way below} $U$, and write
$V\ll U$, if every cover of $U$ by open sets contains a finite subfamily covering $V$. The space $X$
is \emph{core compact} if for every $U\in\OO(X)$ and every $x\in U$ there is $V\in\OO(X)$ with
$x\in V\ll U$. A space is \emph{supercompact} if every open cover of it contains a member that covers
the whole space, and \emph{locally supercompact} if every point has a supercompact open
neighbourhood.
\end{definition}

The relation $\ll$ is transitive and stable under finite unions: if $V_1\ll U$ and $V_2\ll U$, then
$V_1\cup V_2\ll U$, because a finite subcover for $V_1$ together with a finite subcover for $V_2$
covers $V_1\cup V_2$. Core compactness can therefore be read as an approximation property: every
open set is the union of the open sets way below it. It is equivalent to the continuity of the
lattice $\OO(X)$ \citep{hofmann1978spectral,gierz2003continuous,vickers1989topology}, a fact we use in
\S\ref{sec:heredity}.

\begin{definition}\label{def:coherent}
A space $X$ is \emph{coherent} if the intersection of two compact saturated subsets of $X$ is again
compact, and \emph{co-sober} if every nonempty compact saturated set that is not the union of two
proper compact saturated subsets is the saturation of some point \citep{xu2026wilker}.
\end{definition}

Coherence is the order-theoretic counterpart of the requirement that the family of compact
saturated sets be closed under binary intersections, while co-sobriety is the dual of sobriety,
obtained by replacing irreducible closed sets with compact saturated sets that are not unions of two
proper such sets.

\subsection{The lattice of open sets as a topological space}

For $E\subseteq X$ put $\Phi(E)=\{U\in\OO(X): E\subseteq U\}$ and consider the families
\[
  \mathcal{B}_k=\{\Phi(K):K\in\KK(X)\},\quad \mathcal{B}_L=\{\Phi(L):L\in\LL(X)\}.
\]
The members of $\mathcal{B}_k$ and of $\mathcal{B}_L$ are the ``upper'' sets determined by compact,
respectively Lindel\"of, tests. By Lemma~\ref{lem:perm}(i) the union of two compact sets is compact
and the union of two Lindel\"of sets is Lindel\"of, so that both families are closed under finite
intersections; this is the content of the next lemma, which also records the specialization order of
the resulting spaces.

\begin{lemma}\label{lem:bases}
Let $X$ be a space and $E,F\subseteq X$. Then $\Phi(E)\cap\Phi(F)=\Phi(E\cup F)$. Consequently
$\mathcal{B}_k$ and $\mathcal{B}_L$ are bases for topologies on $\OO(X)$; the resulting spaces are
denoted by $\Ok(X)$ and $\OL(X)$. In both of them the specialization order is inclusion, and
$\cl\{U\}=\{V\in\OO(X): V\subseteq U\}$ for every $U\in\OO(X)$.
\end{lemma}

\begin{proof}
For $W\in\OO(X)$ we have $W\in\Phi(E)\cap\Phi(F)$ if and only if $E\subseteq W$ and $F\subseteq W$,
that is, if and only if $E\cup F\subseteq W$; this proves the identity. Since
$\Phi(\emptyset)=\OO(X)$, the members of $\mathcal{B}_k$ cover $\OO(X)$, and the same holds for
$\mathcal{B}_L$. If $K_1,K_2\in\KK(X)$, then
$\Phi(K_1)\cap\Phi(K_2)=\Phi(K_1\cup K_2)\in\mathcal{B}_k$ because $K_1\cup K_2$ is compact, so
$\mathcal{B}_k$ is closed under finite intersections and hence a base; the argument for
$\mathcal{B}_L$ is the same, with Lemma~\ref{lem:perm}(i) replacing compactness.

For the specialization order, let $U,V\in\OO(X)$. A basic neighbourhood of $U$ is a set $\Phi(K)$
with $K\in\KK(X)$ and $K\subseteq U$, and $\Phi(K)$ contains $V$ exactly when $K\subseteq V$. Hence
$U\in\cl\{V\}$ if and only if $K\subseteq V$ for every compact $K\subseteq U$. Taking $K=\{x\}$ for
$x\in U$ shows that $U\subseteq V$ is necessary; it is also sufficient, since $K\subseteq U\subseteq V$
for every $K\subseteq U$. Thus $U\le V$ if and only if $U\subseteq V$, and the computation also
gives $\cl\{U\}=\{V: V\subseteq U\}$.
\end{proof}

The two spaces $\Ok(X)$ and $\OL(X)$ carry the same underlying set $\OO(X)$ but different
topologies, and $\OL(X)$ has the finer one because $\KK(X)\subseteq\LL(X)$. They are the spaces on
which the union operation of the introduction is defined; we write $\mu$ for the binary operation
$\mu(U,V)=U\cup V$ and, for $K\in\KK(X)$, for its restriction
$\mu_K\colon \OL(K)\times\OL(K)\to\Ok(K)$.

%% file: sec3.tex
\section{Wilker-Type Properties and Their Heredity}\label{sec:heredity}

The four properties studied in this paper differ only in the class of sets that is allowed to witness
a decomposition of a compact set along a binary open cover. We first give the unified definition,
then describe the elementary relations between the four properties, and finally study their
behaviour under subspaces, retracts, sums and products.

\begin{definition}\label{def:wilker}
Let $X$ be a space. For $K\in\KK(X)$ and $U_1,U_2\in\OO(X)$ with $K\subseteq U_1\cup U_2$ consider
the condition
\begin{equation}\label{eq:wilker}
  (\exists L_1,L_2\in\mathcal{C})\quad K\subseteq L_1\cup L_2 \ \text{ and }\ L_j\subseteq U_j\
  (j=1,2).
\end{equation}
If \eqref{eq:wilker} holds for all such $K,U_1,U_2$, then, according to the choice of the class
$\mathcal{C}$, the space $X$ is said to have the property listed in Table~\ref{tab:wilker}. The
corresponding classes of spaces are denoted by $\W$, $\Ww$, $\Wi$ and $\Win$.
\end{definition}

\begin{table*}[t]
\centering
\small
\begin{tabular}{@{}lllll@{}}
\toprule
$\mathcal{C}$ & $\KK(X)$ & $\LL(X)$ & $\LL(K)$ & $\KK(K)$\\
\midrule
Property & Wilker $\W$ & weak Wilker $\Ww$ & inner weak Wilker $\Wi$ & inner compact refinement $\Win$\\
\bottomrule
\end{tabular}
\caption{The four Wilker-type properties, according to the class $\mathcal{C}$ of admissible
witnesses in condition~\eqref{eq:wilker}.}
\label{tab:wilker}
\end{table*}

When $\mathcal{C}=\LL(K)$ or $\mathcal{C}=\KK(K)$, the witnesses must lie inside $K$, and then the
inclusion $K\subseteq L_1\cup L_2$ of \eqref{eq:wilker} is equivalent to the identity
$K=L_1\cup L_2$. This is the reason for the word ``inner'': the decomposition is a decomposition of
$K$ itself and not only of a neighbourhood of $K$.

\begin{remark}\label{rem:hierarchy}
Since every compact set is Lindel\"of, $\W\subseteq\Ww$, and since $\LL(K)\subseteq\LL(X)$ and
$\KK(K)\subseteq\LL(K)$, we have
\[
  \Win\ \subseteq\ \W\cap\Wi\ \subseteq\ \Ww .
\]
The properties $\W$ and $\Wi$ are not comparable: Proposition~\ref{prop:M0} below produces a space in
$\Wi\setminus\W$, and Proposition~\ref{prop:locallysuper} produces a space in $\W\setminus\Wi$. Hence
the only inclusions among the four classes are the two displayed above, together with the inclusion
$\Win\subseteq\Ww$ that follows from them.
\end{remark}

We now show that $\Wi$ can be tested on compact subspaces only. This is the technical reason why the
inner weak Wilker property is much better behaved than $\Ww$.

\begin{proposition}\label{prop:Wi-char}
For a space $X$ the following are equivalent:
\begin{enumerate}
\item $X\in\Wi$;
\item every subspace $B\subseteq X$ belongs to $\Wi$;
\item every subspace $B\subseteq X$ belongs to $\Ww$;
\item every compact subspace $K\in\KK(X)$ belongs to $\Ww$.
\end{enumerate}
\end{proposition}

\begin{proof}
$(1)\Rightarrow(2)$. Assume $X\in\Wi$ and let $B\subseteq X$ be a subspace. Take $K\in\KK(B)$ and
open sets $U_1,U_2\in\OO(X)$ with $K\subseteq (B\cap U_1)\cup (B\cap U_2)$; such a pair exists for
every binary open cover of $K$ in $B$, since each open set of $B$ is of the form $B\cap U$. Then $K$
is compact in $X$ by Lemma~\ref{lem:perm}(ii), and $K\subseteq U_1\cup U_2$. So $\Wi$ supplies
$L_1,L_2\in\LL(X)$ with $L_j\subseteq K\cap U_j$ and $K=L_1\cup L_2$. From $L_j\subseteq K\subseteq B$
we obtain $L_j\subseteq B\cap U_j$, and $L_j$ is Lindel\"of as a subspace of $B$ by
Lemma~\ref{lem:perm}(ii). Hence $B\in\Wi$.

$(2)\Rightarrow(3)$. This is immediate from $\Wi\subseteq\Ww$.

$(3)\Rightarrow(4)$. A compact subspace is a subspace, so it belongs to $\Ww$.

$(4)\Rightarrow(1)$. Assume that every compact subspace of $X$ belongs to $\Ww$, and take
$K\in\KK(X)$, $U_1,U_2\in\OO(X)$ with $K\subseteq U_1\cup U_2$. The sets $K\cap U_1$ and $K\cap U_2$
are open in the subspace $K$ and cover it, so, since $K\in\Ww$, there are $L_1,L_2\in\LL(K)$ with
$L_j\subseteq K\cap U_j$ and $K=L_1\cup L_2$. By Lemma~\ref{lem:perm}(ii) each $L_j$ is Lindel\"of as
a subspace of $X$, and $L_j\subseteq U_j$. Hence $X\in\Wi$.
\end{proof}

\begin{proposition}\label{prop:heredity}
The properties $\Wi$ and $\Win$ are hereditary with respect to arbitrary subspaces.
\end{proposition}

\begin{proof}
The statement for $\Wi$ is part of Proposition~\ref{prop:Wi-char}. For $\Win$, let $X\in\Win$, let
$B\subseteq X$ and let $K\in\KK(B)$, $U_1,U_2\in\OO(X)$ with
$K\subseteq (B\cap U_1)\cup (B\cap U_2)$. Then $K\in\KK(X)$ and $K\subseteq U_1\cup U_2$, so $\Win$
gives $F_1,F_2\in\KK(X)$ with $F_j\subseteq K\cap U_j$ and $K=F_1\cup F_2$. Since
$F_j\subseteq K\subseteq B$, each $F_j$ is compact in the subspace $B$ by Lemma~\ref{lem:perm}(ii),
and $F_j\subseteq B\cap U_j$. Hence $B\in\Win$.
\end{proof}

The corresponding statement for $\W$ and $\Ww$ is false for arbitrary subspaces, but it holds for
closed ones.

\begin{proposition}\label{prop:closed}
The properties $\W$ and $\Ww$ are hereditary with respect to closed subspaces.
\end{proposition}

\begin{proof}
Let $B$ be a closed subspace of $X$, let $X\in\Ww$ and let $K\in\KK(B)$, $U_1,U_2\in\OO(X)$ with
$K\subseteq (B\cap U_1)\cup (B\cap U_2)$. Then $K\subseteq U_1\cup U_2$, so there are
$L_1,L_2\in\LL(X)$ with $L_j\subseteq U_j$ and $K\subseteq L_1\cup L_2$. Put $L'_j=L_j\cap B$. Since
$B$ is closed in $X$, the set $L'_j$ is closed in the Lindel\"of space $L_j$, hence Lindel\"of by
Lemma~\ref{lem:perm}(iii), and $L'_j\subseteq B\cap U_j$. Moreover $K\subseteq B$ gives
\[
  K=K\cap (L_1\cup L_2)\subseteq (L_1\cap B)\cup (L_2\cap B)=L'_1\cup L'_2 .
\]
Therefore $B\in\Ww$. The argument for $\W$ is the same, with ``compact'' in place of ``Lindel\"of'':
a closed subset of a compact space is compact.
\end{proof}

That the restriction to closed subspaces in Proposition~\ref{prop:closed} cannot be removed is shown
by Proposition~\ref{prop:locallysuper}, where a locally supercompact space in $\W\setminus\Wi$ is
constructed; its subspace $M$ belongs to $\W$ but not to $\Ww$, although the ambient space does.

\begin{proposition}\label{prop:retract}
The properties $\W$, $\Ww$ and $\Wi$ are preserved by continuous retracts: if $e\colon X\to Y$ is an
embedding and $r\colon Y\to X$ is a continuous map with $r\circ e=\id_X$, then $Y\in\mathbb{W}$
implies $X\in\mathbb{W}$ for each $\mathbb{W}\in\{\W,\Ww,\Wi\}$.
\end{proposition}

\begin{proof}
Let $e\colon X\to Y$ be an embedding with retraction $r$, so that $r(e(x))=x$ for all $x\in X$.

Assume first that $Y\in\Ww$, and take $K\in\KK(X)$, $U_1,U_2\in\OO(X)$ with $K\subseteq U_1\cup U_2$.
Since $e$ is continuous, $e(K)$ is compact in $Y$; since $r$ is continuous, the sets
$W_j=r^{-1}(U_j)$ are open in $Y$. For $k\in K$ we have $r(e(k))=k\in U_1\cup U_2$, hence
$e(k)\in W_1\cup W_2$; therefore $e(K)\subseteq W_1\cup W_2$. As $Y\in\Ww$, there are
$L_1,L_2\in\LL(Y)$ with $L_j\subseteq W_j$ and $e(K)\subseteq L_1\cup L_2$. Put $L'_j=r(L_j)$. Then
$L'_j$ is Lindel\"of, being a continuous image of $L_j$, and $L_j\subseteq r^{-1}(U_j)$ gives
$L'_j\subseteq U_j$. Finally
\[
  K=r(e(K))\subseteq r(L_1\cup L_2)=L'_1\cup L'_2 .
\]
Hence $X\in\Ww$.

Assume next that $Y\in\Wi$. With $K,U_1,U_2,W_1,W_2$ as above we have $e(K)\subseteq W_1\cup W_2$, so
there are $L_1,L_2\in\LL(Y)$ with $L_j\subseteq e(K)\cap W_j$ and $e(K)=L_1\cup L_2$. Putting
$L'_j=r(L_j)$ we obtain, as before, that $L'_j$ is Lindel\"of, and now
\[
  L'_j\subseteq r\bigl(e(K)\cap r^{-1}(U_j)\bigr)\subseteq r(e(K))\cap U_j=K\cap U_j ,
\]
while $K=r(e(K))=L'_1\cup L'_2$. Hence $X\in\Wi$.

The case of $\W$ is obtained from the first part by reading ``compact'' for ``Lindel\"of''
everywhere: the continuous image of a compact set is compact.
\end{proof}

\begin{proposition}\label{prop:sums}
Let $\mathbb{W}\in\{\Ww,\Wi,\Win\}$. For a family $\{X_\alpha\}_{\alpha\in I}$ of nonempty spaces,
the topological sum $\bigsqcup_{\alpha\in I}X_\alpha$ belongs to $\mathbb{W}$ if and only if
$X_\alpha\in\mathbb{W}$ for every $\alpha\in I$.
\end{proposition}

\begin{proof}
\emph{Necessity.} If $X=\bigsqcup_{\alpha}X_\alpha\in\mathbb{W}$, then each summand $X_\alpha$ is a
clopen subspace of $X$. For $\mathbb{W}=\Ww$ we apply Proposition~\ref{prop:closed} (closed
subspaces), and for $\mathbb{W}=\Wi$ or $\mathbb{W}=\Win$ we apply
Proposition~\ref{prop:heredity} (arbitrary subspaces); in each case $X_\alpha\in\mathbb{W}$.

\emph{Sufficiency.} Assume $X_\alpha\in\mathbb{W}$ for all $\alpha$, and let
$X=\bigsqcup_{\alpha}X_\alpha$. Let $K\in\KK(X)$ and let $U_1,U_2\in\OO(X)$ with
$K\subseteq U_1\cup U_2$. Since $\{X_\alpha\}$ is a pairwise disjoint open cover of $X$, the compact
set $K$ meets only finitely many summands; write $I_K=\{\alpha_1,\dots,\alpha_n\}$ for the finite set
of indices with $K\cap X_\alpha\neq\emptyset$. For each $\alpha\in I_K$ the set $K\cap X_\alpha$ is
compact in the subspace $X_\alpha$, and
\[
  K\cap X_\alpha\subseteq (X_\alpha\cap U_1)\cup (X_\alpha\cap U_2).
\]
Since $X_\alpha\in\mathbb{W}$, there are sets $L_{1,\alpha},L_{2,\alpha}$ in $X_\alpha$, of the kind
required by $\mathbb{W}$, with $L_{j,\alpha}\subseteq X_\alpha\cap U_j$ and
$K\cap X_\alpha\subseteq L_{1,\alpha}\cup L_{2,\alpha}$; when $\mathbb{W}\in\{\Wi,\Win\}$ we may and
do choose them inside $K\cap X_\alpha$.

If $\mathbb{W}=\Ww$, put $L_j=\bigcup_{\alpha\in I_K}L_{j,\alpha}$. This is a Lindel\"of subset of
$X$ by Lemma~\ref{lem:perm}(i), because $I_K$ is finite, and $L_j\subseteq U_j$; moreover
$K\subseteq L_1\cup L_2$. Hence $X\in\Ww$. If $\mathbb{W}=\Wi$, the sets $L_{j,\alpha}$ are
Lindel\"of and $L_j=\bigcup_{\alpha\in I_K}L_{j,\alpha}\subseteq K\cap U_j$ is Lindel\"of as a finite
union of Lindel\"of sets, so $X\in\Wi$. If $\mathbb{W}=\Win$, the sets $L_{j,\alpha}$ are compact
subsets of $K\cap U_j$ and their finite union is compact, so $X\in\Win$.
\end{proof}

The next two results are the technical core of the section: they show that over KC-spaces the weak
and the inner weak Wilker properties coincide, and that the inner compact refinement property is
automatic as soon as all compact subspaces are Hausdorff.

\begin{proposition}\label{prop:KC}
Let $X$ be a space.
\begin{enumerate}
\item If $X$ is a KC-space, then $X\in\Ww$ if and only if $X\in\Wi$.
\item If every compact subspace of $X$ is Hausdorff, then $X\in\Win$.
\end{enumerate}
\end{proposition}

\begin{proof}
(i) Let $X$ be a KC-space with $X\in\Ww$, and take $K\in\KK(X)$, $U_1,U_2\in\OO(X)$ with
$K\subseteq U_1\cup U_2$. By $\Ww$ there are $L_1,L_2\in\LL(X)$ with $L_j\subseteq U_j$ and
$K\subseteq L_1\cup L_2$. Since $X$ is KC, the compact set $K$ is closed in $X$, so $L_j\cap K$ is a
closed subset of the Lindel\"of space $L_j$ and hence Lindel\"of. Moreover
$L_j\cap K\subseteq K\cap U_j$, and
\[
  K=K\cap (L_1\cup L_2)=(L_1\cap K)\cup (L_2\cap K).
\]
Thus the sets $L'_j=L_j\cap K$ are Lindel\"of subsets of $K$ contained in $K\cap U_j$ and they cover
$K$, which is exactly $X\in\Wi$. The converse is Remark~\ref{rem:hierarchy}.

(ii) Assume that every compact subspace of $X$ is Hausdorff, and take $K\in\KK(X)$,
$U_1,U_2\in\OO(X)$ with $K\subseteq U_1\cup U_2$. Put $C_j=K\setminus U_j$; these are closed subsets
of the compact space $K$, and $C_1\cap C_2=\emptyset$ because $U_1\cup U_2$ covers $K$. By
hypothesis $K$ is compact Hausdorff, hence normal, so there are relatively open sets $G_1,G_2$ of
$K$ with $C_j\subseteq G_j$ and $G_1\cap G_2=\emptyset$. Let $F_j=K\setminus G_j$. Each $F_j$ is
closed in the compact space $K$ and therefore compact, and
\[
  \begin{gathered}
    F_j\subseteq K\setminus C_j=K\cap U_j,\\
    F_1\cup F_2=(K\setminus G_1)\cup (K\setminus G_2)=K\setminus (G_1\cap G_2)=K .
  \end{gathered}
\]
Hence $X\in\Win$.
\end{proof}

We now separate $\Wi$ from $\W$. The examples are countably infinite compact KC-spaces; in a
countable space every subset is Lindel\"of, so that $\Wi$ holds automatically, and the KC property
together with hyperconnectedness prevents a decomposition into two proper closed pieces, which
excludes $\W$.

\begin{proposition}\label{prop:M0}
There is a countably infinite compact KC-space $M_0$ with $M_0\in\Wi\setminus\W$.
\end{proposition}

\begin{proof}
Let $M_0$ be Kov\'ar's countably infinite compact Fr\'echet--Urysohn hyperconnected KC-space
\citep[Example~2.2 and the following paragraph, p.~6]{kovar2005maximal}; the construction is recalled
in~\cite[before Lemma~3.14]{xu2026wilker}.

\emph{$M_0\in\Wi$.} Since $M_0$ is countable, every subset of it is countable, and every countable
space is Lindel\"of. Given $K\in\KK(M_0)$ and $U_1,U_2\in\OO(M_0)$ with $K\subseteq U_1\cup U_2$, put
$L_j=K\cap U_j$. Each $L_j$ is a subset of $K$, hence a Lindel\"of subset of $M_0$, and
$L_j\subseteq K\cap U_j$; moreover
\[
  L_1\cup L_2=(K\cap U_1)\cup (K\cap U_2)=K\cap (U_1\cup U_2)=K .
\]
So the witnesses required by $\Wi$ exist.

\emph{$M_0\notin\W$.} Suppose that $M_0\in\W$. Choose distinct points $a,b\in M_0$ and put
$U_1=M_0\setminus\{a\}$, $U_2=M_0\setminus\{b\}$. Both are open: for instance $\{a\}$ is compact and
hence closed, because $M_0$ is a KC-space, so its complement $U_1$ is open. Clearly
$U_1\cup U_2=M_0$, so the compact set $K=M_0$ satisfies $K\subseteq U_1\cup U_2$. By $\W$ there are
compact sets $K_1\subseteq U_1$ and $K_2\subseteq U_2$ with $M_0=K_1\cup K_2$. Being compact in the
KC-space $M_0$, the sets $K_1$ and $K_2$ are closed; and $K_1\subseteq M_0\setminus\{a\}$ gives
$a\notin K_1$, so $K_1\neq M_0$, and similarly $K_2\neq M_0$. Hence $M_0\setminus K_1$ and
$M_0\setminus K_2$ are nonempty open sets, while
\[
  (M_0\setminus K_1)\cap (M_0\setminus K_2)=M_0\setminus (K_1\cup K_2)=\emptyset,
\]
contradicting the hyperconnectedness of $M_0$. Therefore $M_0\notin\W$, and $M_0\in\Wi\setminus\W$.
\end{proof}

The space below is a more concrete example with the same effect; it is obtained from $\QQ$ by adding
a point at infinity whose neighbourhoods are the complements of compact subsets of $\QQ$.

\begin{example}\label{ex:Qplus}
Let $\QQ^{+}=\QQ\cup\{\infty\}$, where $\QQ$ carries its usual topology and the neighbourhoods of
$\infty$ are the sets $\QQ^{+}\setminus C$ with $C\in\KK(\QQ)$. Then $\QQ^{+}$ is a countably
infinite compact KC-space which is not Hausdorff, and $\QQ^{+}\in\Wi\setminus\W$.
\end{example}

\begin{proof}
The family $\{\QQ^{+}\setminus C:C\in\KK(\QQ)\}$ is a neighbourhood base at $\infty$: it is closed
under finite intersections because the union of two compact subsets of $\QQ$ is compact, and each of
its members contains $\infty$.

\emph{$\QQ^{+}$ is compact.} Let $\mathcal{U}$ be an open cover of $\QQ^{+}$ and let
$U\in\mathcal{U}$ contain $\infty$. Then $U$ contains some $\QQ^{+}\setminus C$ with $C\in\KK(\QQ)$,
so $C$ is covered by the remaining members of $\mathcal{U}$; since $C$ is compact in $\QQ$, finitely
many of them cover $C$. Together with $U$ this gives a finite subcover.

\emph{$\QQ^{+}$ is a KC-space.} For $C\in\KK(\QQ)$ the set $\QQ^{+}\setminus C$ is open, so $C$ is
closed in $\QQ^{+}$. Now let $E$ be a compact subset of $\QQ^{+}$. If $\infty\notin E$, then
$E\subseteq \QQ$, and $E$ is compact in the open subspace $\QQ$ of $\QQ^{+}$, hence compact in $\QQ$
and therefore closed in $\QQ$; since $\QQ$ is open in $\QQ^{+}$, it follows that $E$ is closed in
$\QQ^{+}$. If $\infty\in E$, take $x\in\cl_{\QQ}(E\cap\QQ)$. As $\QQ$ is first countable, there are
$x_n\in E\cap\QQ$ with $x_n\to x$. The set $H=\{x\}\cup\{x_n: n<\omega\}$ is compact in $\QQ$, being a
convergent sequence together with its limit, so $H$ is closed in $\QQ^{+}$. Hence $E\cap H$ is a
closed subset of the compact set $E$ and is therefore compact; it is contained in the Hausdorff space
$H$, so it is closed in $H$. From $x_n\in E\cap H$ and $x_n\to x$ we obtain
$x\in\cl_H(E\cap H)=E\cap H\subseteq E$. Consequently $E\cap\QQ$ is closed in $\QQ$, and $E$ is
closed in $\QQ^{+}$. Thus $\QQ^{+}$ is a KC-space, and it is countable by construction.

\emph{$\QQ^{+}$ is not Hausdorff.} Suppose that $0\in U$ and $\infty\in V$ are disjoint open sets.
Since $V$ contains some $\QQ^{+}\setminus C$ with $C\in\KK(\QQ)$, we obtain $U\subseteq C$. On the other
hand $U$ contains a rational neighbourhood of $0$, so there is $\varepsilon>0$ with
$(-\varepsilon,\varepsilon)\cap\QQ\subseteq U\subseteq C$. The set $C$ is compact as a subspace of
$\QQ$, hence compact and closed in $\RR$, so
$[-\varepsilon,\varepsilon]=\cl_{\RR}\bigl((-\varepsilon,\varepsilon)\cap\QQ\bigr)\subseteq C\subseteq \QQ$,
which is absurd because $[-\varepsilon,\varepsilon]$ contains irrational numbers.

\emph{$\QQ^{+}\in\Wi$.} Since $\QQ^{+}$ is countable, every subset of it is Lindel\"of, and the
argument used for $M_0$ in the proof of Proposition~\ref{prop:M0} applies verbatim: for
$K\in\KK(\QQ^{+})$ and $K\subseteq U_1\cup U_2$, the sets $L_j=K\cap U_j$ are Lindel\"of and satisfy
$L_1\cup L_2=K$.

\emph{$\QQ^{+}\notin\W$.} Suppose that $\QQ^{+}\in\W$, and put $U_1=\QQ^{+}\setminus\{0\}$ and
$U_2=\QQ^{+}\setminus\{\infty\}=\QQ$. Both are open: the singleton $\{0\}$ is compact and hence
closed, because $\QQ^{+}$ is a KC-space. Since $U_1\cup U_2=\QQ^{+}$, the compact set $K=\QQ^{+}$
satisfies $K\subseteq U_1\cup U_2$, so by $\W$ there are compact sets $K_1\subseteq U_1$ and
$K_2\subseteq U_2$ with $\QQ^{+}=K_1\cup K_2$. Being compact in a KC-space, they are closed, and
$0\notin K_1$, $\infty\notin K_2$, so both are proper closed subsets. Hence
$\QQ^{+}\setminus K_1$ and $\QQ^{+}\setminus K_2$ are nonempty open sets containing $0$ and $\infty$
respectively, and their intersection is empty, a contradiction. Therefore $\QQ^{+}\notin\W$, and
$\QQ^{+}\in\Wi\setminus\W$.
\end{proof}

The following theorem identifies $\Ww$ and $\Wi$ with continuity conditions for the union operation
on the lattices of open sets; it is the reason why these properties are relevant for power spaces.
The proof uses the fact, recorded in Lemma~\ref{lem:perm}(v), that a basic neighbourhood of an open
set may always be taken to be generated by a saturated Lindel\"of set.

\begin{theorem}\label{thm:mu}
Let $X$ be a space. Then $X\in\Ww$ if and only if $\mu\colon \OL(X)\times\OL(X)\to\Ok(X)$,
$\mu(U,V)=U\cup V$, is continuous. Moreover $X\in\Wi$ if and only if, for every $K\in\KK(X)$, the
restriction $\mu_K\colon \OL(K)\times\OL(K)\to\Ok(K)$, $\mu_K(U,V)=U\cup V$, is continuous.
\end{theorem}

\begin{proof}
We use the bases of Lemma~\ref{lem:bases}. Observe first that $\Phi(L)=\Phi(\sat L)$ and that
$\sat L$ is Lindel\"of whenever $L$ is, by Lemma~\ref{lem:perm}(v); hence a basic open set of
$\Ok(X)$ or of $\OL(X)$ can always be written as $\Phi(L)$ with $L$ Lindel\"of and saturated.

\emph{The case of $\Ww$.} Assume $X\in\Ww$ and let $\Phi(K)$ be a basic open set of $\Ok(X)$, with
$K\in\KK(X)$. Then $\mu^{-1}\Phi(K)=\{(U,V):K\subseteq U\cup V\}$. Given $(U,V)$ in this set, $\Ww$
provides $L_1,L_2\in\LL(X)$ with $L_j\subseteq U_j$ (where $U_1=U$, $U_2=V$) and
$K\subseteq L_1\cup L_2$. Then $(U,V)\in\Phi(L_1)\times\Phi(L_2)$, and every
$(U',V')\in\Phi(L_1)\times\Phi(L_2)$ satisfies $K\subseteq L_1\cup L_2\subseteq U'\cup V'$. Hence
$\Phi(L_1)\times\Phi(L_2)\subseteq\mu^{-1}\Phi(K)$, and $\mu$ is continuous at $(U,V)$.

Conversely, assume that $\mu$ is continuous, and take $K\in\KK(X)$ and $U_1,U_2\in\OO(X)$ with
$K\subseteq U_1\cup U_2$. Then $(U_1,U_2)\in\mu^{-1}\Phi(K)$, so by continuity there are
$L_1,L_2\in\LL(X)$ with $(U_1,U_2)\in\Phi(L_1)\times\Phi(L_2)$ and
$\Phi(L_1)\times\Phi(L_2)\subseteq\mu^{-1}\Phi(K)$. In particular $L_j\subseteq U_j$, and every pair
of open sets $G\supseteq L_1$, $H\supseteq L_2$ satisfies $K\subseteq G\cup H$. We claim that
$K\subseteq \sat L_1\cup \sat L_2$. Otherwise there is $x\in K\setminus(\sat L_1\cup\sat L_2)$; since
$x\notin\sat L_1=\bigcap\{G\in\OO(X):L_1\subseteq G\}$, there is an open set $G\supseteq L_1$ with
$x\notin G$, and similarly there is an open set $H\supseteq L_2$ with $x\notin H$. Then
$K\subseteq G\cup H$ but $x\notin G\cup H$, a contradiction. Since $\sat L_j\subseteq U_j$ and
$\sat L_j$ is Lindel\"of, we conclude that $X\in\Ww$.

\emph{The case of $\Wi$.} Assume $X\in\Wi$, let $K\in\KK(X)$ and let $\Phi_K(F)$ be a basic open set
of $\Ok(K)$, where $F\in\KK(K)$ and $\Phi_K(E)=\{W\in\OO(K):E\subseteq W\}$ for $E\subseteq K$. Take
$(U,V)\in\mu_K^{-1}\Phi_K(F)$, so that $F\subseteq U\cup V$. Choose $U',V'\in\OO(X)$ with
$U=K\cap U'$ and $V=K\cap V'$; then $F\subseteq U'\cup V'$, and since $F\in\KK(X)$ we may apply
$X\in\Wi$ and obtain $L_1,L_2\in\LL(X)$ with $L_j\subseteq F\cap U'_j$ and $F=L_1\cup L_2$. Then
$L_j\subseteq F\subseteq K$, so $L_j\in\LL(K)$ and $L_j\subseteq U'_j\cap K=U_j$. Thus
$(U,V)\in\Phi_K(L_1)\times\Phi_K(L_2)$, and every $(U'',V'')\in\Phi_K(L_1)\times\Phi_K(L_2)$
satisfies $F=L_1\cup L_2\subseteq U''\cup V''$, so that
$\Phi_K(L_1)\times\Phi_K(L_2)\subseteq\mu_K^{-1}\Phi_K(F)$. Hence $\mu_K$ is continuous.

Conversely, assume that $\mu_K$ is continuous for every $K\in\KK(X)$, and take $K\in\KK(X)$,
$U_1,U_2\in\OO(X)$ with $K\subseteq U_1\cup U_2$. Put $V_j=K\cap U_j$; then $V_1,V_2$ are open in $K$
and $K=V_1\cup V_2$, so $(V_1,V_2)\in\mu_K^{-1}\Phi_K(K)$. By continuity there are $L_1,L_2\in\LL(K)$
with $L_j\subseteq V_j$ and $\Phi_K(L_1)\times\Phi_K(L_2)\subseteq\mu_K^{-1}\Phi_K(K)$. This means
that $K\subseteq W_1\cup W_2$ whenever $W_j$ are open in $K$ with $W_j\supseteq L_j$. Repeating the
argument of the previous paragraph inside the space $K$, with
$\sat_K L=\bigcap\{W\in\OO(K):L\subseteq W\}$ in place of $\sat L$, we obtain
$K\subseteq\sat_K L_1\cup\sat_K L_2$. Now $\sat_K L_j\subseteq V_j$ and $\sat_K L_j$ is Lindel\"of as
a subspace of $K$, hence as a subspace of $X$; since $V_j\subseteq U_j$ and
$K\subseteq\sat_K L_1\cup\sat_K L_2$, we conclude that $X\in\Wi$.
\end{proof}

Core compactness gives a sufficient condition for the weak Wilker property. The proof rests on the
following standard interpolation property of the way-below relation; it holds because a core compact
space has a continuous lattice of open sets \citep{gierz2003continuous,vickers1989topology}.

\begin{lemma}\label{lem:interp}
Let $X$ be a core compact space and let $V\ll U$ be open subsets of $X$. Then there is an open set $W$
with $V\ll W\ll U$.
\end{lemma}

\begin{proof}
For a core compact space the lattice $\OO(X)$, ordered by inclusion, is continuous
\citep{gierz2003continuous}, and in a continuous lattice the way-below relation is interpolative
\citep[Lemma~I-1.20]{gierz2003continuous}: from $V\ll U$ one obtains an open set $W$ with
$V\ll W\ll U$. Note that $W\subseteq U$, since $W\ll U$ forces $W\subseteq U$.
\end{proof}

\begin{theorem}\label{thm:corecompact}
Let $X$ be a core compact space. Then every open subset of $X$ is a union of Lindel\"of open sets;
more precisely, for every open cover $\mathcal{U}$ of a compact set $K\subseteq X$ there are finitely
many Lindel\"of open sets, each contained in a member of $\mathcal{U}$, that cover $K$. In particular
$X\in\Ww$.
\end{theorem}

\begin{proof}
\emph{Step 1: a point has arbitrarily small Lindel\"of open neighbourhoods.} Fix $x\in U$ with $U$
open. By core compactness there is $V_0\in\OO(X)$ with $x\in V_0\ll U$. Applying
Lemma~\ref{lem:interp} repeatedly, we obtain open sets $V_0,V_1,V_2,\dots$ with
\[
  x\in V_0\ll V_1\ll V_2\ll\cdots\ll U .
\]
Put $N=\bigcup_{n<\omega}V_n$. Then $N$ is open and $x\in N\subseteq U$. To see that $N$ is
Lindel\"of, let $\mathcal{G}$ be a cover of $N$ by open sets. For fixed $n$, the family of finite
unions of members of $\mathcal{G}$ is directed and its union contains $N$, hence contains $V_{n+1}$;
since $V_n\ll V_{n+1}$, there is a finite $\mathcal{G}_n\subseteq\mathcal{G}$ with
$V_n\subseteq\bigcup\mathcal{G}_n$. Therefore
$N=\bigcup_n V_n\subseteq\bigcup_n\bigcup\mathcal{G}_n$, and $\bigcup_n\mathcal{G}_n$ is a countable
subfamily of $\mathcal{G}$ covering $N$. Thus $N$ is a Lindel\"of open set with $x\in N\subseteq U$.

\emph{Step 2: finitely many such sets cover a compact set.} Let $\mathcal{U}$ be an open cover of
$K\in\KK(X)$. For each $x\in K$ choose $U_x\in\mathcal{U}$ with $x\in U_x$ and, by Step~1, a
Lindel\"of open set $N_x$ with $x\in N_x\subseteq U_x$. Since $K$ is compact, there is a finite
$F\subseteq K$ with $K\subseteq\bigcup_{x\in F}N_x$. Grouping the finitely many sets $N_x$ according
to the member $U_x$ of $\mathcal{U}$ that contains them, and taking finite unions, we obtain finitely
many Lindel\"of open sets, each contained in a member of $\mathcal{U}$, whose union contains $K$.

\emph{Step 3: $X\in\Ww$.} Given $K\in\KK(X)$ and $K\subseteq U_1\cup U_2$, apply Step~2 to the
two-element cover $\{U_1,U_2\}$: there are Lindel\"of open sets $L_1\subseteq U_1$,
$L_2\subseteq U_2$ with $K\subseteq L_1\cup L_2$. Hence $X\in\Ww$.
\end{proof}

The last result of this section describes the behaviour of $\Ww$ under products. Its hypothesis is
that one factor belongs to $\Ww$ and that the compact subspaces of the other factor are Hausdorff;
under this assumption the product remains in $\Ww$, and if moreover the first factor is KC, then the
product is KC and belongs to $\Wi$.

\begin{theorem}\label{thm:product}
Let $X\in\Ww$ and suppose that every compact subspace of $Y$ is Hausdorff. Then $X\times Y\in\Ww$. If
in addition $X$ is a KC-space, then $X\times Y$ is a KC-space and $X\times Y\in\Wi$.
\end{theorem}

\begin{proof}
Let $K\in\KK(X\times Y)$ and let $U_1,U_2$ be open subsets of $X\times Y$ with
$K\subseteq U_1\cup U_2$; we may assume $K\neq\emptyset$. Put $C=\pi_Y(K)$, a compact subset of $Y$,
and replace each $U_i$ by $U_i\cap (X\times C)$, which is open in $X\times C$ and which we still
denote by $U_i$.

\emph{Fibres over $C$.} Since every compact subspace of $Y$ is Hausdorff, $C$ is Hausdorff, hence so
is $X\times C$. For $y\in C$ the set $K_y=K\cap (X\times\{y\})$ is closed in $K$, because
$X\times\{y\}$ is closed in $X\times C$, and is therefore compact. As $K_y\subseteq U_1\cup U_2$,
finitely many open rectangles $G_r\times H_r$ with $G_r\times H_r\subseteq U_{c(r)}$,
$c(r)\in\{1,2\}$, cover $K_y$, and we may assume $y\in H_r$ for all $r$. Set
\[
  V_i=\bigcup_{c(r)=i}G_r,\qquad H=\bigcap_r H_r .
\]
Then $y\in H$, $K_y\subseteq (V_1\cup V_2)\times H$, and $V_i\times H\subseteq U_i$ for $i=1,2$.

\emph{Shrinking the base.} Put $D_y=\pi_Y\bigl(K\setminus ((V_1\cup V_2)\times C)\bigr)$. The set
$D_y$ is compact, being a continuous image of a closed subset of the compact space $K$, and
$y\notin D_y$ because $K_y\subseteq (V_1\cup V_2)\times C$. As $C$ is compact Hausdorff, it is
regular, so there is an open set $N_y$ of $C$ with $y\in N_y$ and
$B_y=\cl_C N_y\subseteq C\setminus D_y$. Observe that $B_y\subseteq H$. The set
$T_y=K\cap (X\times B_y)$ is compact, being closed in $K$, and $\pi_X(T_y)\subseteq V_1\cup V_2$ by
the choice of $B_y$; moreover $V_i\times B_y\subseteq U_i$.

\emph{Conclusion.} Since $\pi_X(T_y)$ is compact and $X\in\Ww$, there are $A_1(y),A_2(y)\in\LL(X)$
with $\pi_X(T_y)\subseteq A_1(y)\cup A_2(y)$ and $A_i(y)\subseteq V_i$. The family $\{N_y:y\in C\}$
covers the compact set $C$, so there is a finite $J\subseteq C$ with $C=\bigcup_{y\in J}N_y$. Put
\[
  L_i=\bigcup_{y\in J}A_i(y)\times B_y .
\]
Then $L_i\subseteq U_i$, because $A_i(y)\subseteq V_i$ and $V_i\times B_y\subseteq U_i$, and
$K=\bigcup_{y\in J}T_y\subseteq L_1\cup L_2$. It remains to note that $L_i$ is Lindel\"of: the
product of a Lindel\"of space and a compact space is Lindel\"of, since for such a product $A\times B$
one covers each slice $\{a\}\times B$ by finitely many rectangles, intersects the corresponding
first-coordinate neighbourhoods to obtain an open $G_a\ni a$, chooses a countable $A_0\subseteq A$
with $A=\bigcup_{a\in A_0}G_a$ and collects the countably many finite subfamilies so obtained.
Therefore $L_i$ is a finite union of Lindel\"of sets, hence Lindel\"of, and $X\times Y\in\Ww$.

\emph{The KC case.} Assume now that $X$ is a KC-space. We first note that every compact subset of $Y$
is closed in $Y$: if $D\subseteq Y$ is compact and $z\in\cl_Y D$, then $D\cup\{z\}$ is compact, being
a subset of the compact set $\cl_Y D$, and it is Hausdorff by hypothesis, so $D$ is closed in
$D\cup\{z\}$ and $z\in D$. In particular $C$ is closed in $Y$, and $K$ is a closed subset of
$X\times C$.

To see that $K$ is closed in $X\times C$, let $(x,y)\notin K$. Since $X$ is KC, it is $T_1$, so
$\{x\}$ is closed in $X$ and $K_x=\pi_Y(K\cap(\{x\}\times C))$ is compact; it is closed in the
compact Hausdorff space $C$. Choose a closed neighbourhood $B$ of $y$ in $C$ with
$B\cap K_x=\emptyset$, and put $D=\pi_X(K\cap (X\times B))$. Then $D$ is compact, and it is closed in
$X$ because $X$ is KC; moreover $x\notin D$. Hence
\[
  (x,y)\in (X\setminus D)\times \operatorname{int}_C B\subseteq (X\times C)\setminus K ,
\]
which shows that $K$ is closed. Thus $X\times Y$ is a KC-space, and Proposition~\ref{prop:KC}(i) gives
$X\times Y\in\Wi$.
\end{proof}

%% file: sec4.tex
\section{Counterexamples in Power Spaces and Function Spaces}\label{sec:counter}

The heredity results of \S\ref{sec:heredity} show that the inner properties are stable, but they say
nothing about products and about the hyperspace constructions. In this section we show that all of
these can fail, and we use the failures to answer a question raised in~\cite{xu2026wilker} and to
produce a non-sober function space.

The basic building block is a compact $T_1$ space whose topology is so large that no two Lindel\"of
pieces can cover a binary open cover.

\begin{lemma}\label{lem:grid}
Let $A$ be a set with $|A|=\aleph_2$, put $D=A\times A$ and $M=D\sqcup\{p,q\}$. Declare every point
of $D$ isolated, and take at $p$ and $q$ the neighbourhood bases
\[
  N_p(F)=\{p\}\cup\bigl((A\setminus F)\times A\bigr),\quad
  N_q(G)=\{q\}\cup\bigl(A\times (A\setminus G)\bigr),
\]
where $F,G$ range over the finite subsets of $A$. Then $M$ is compact and $T_1$, and for
$U=M\setminus\{p\}$, $V=M\setminus\{q\}$ there are no $L_U\in\LL(U)$ and $L_V\in\LL(V)$ with
$M\subseteq L_U\cup L_V$. In particular $M\notin\Ww$.
\end{lemma}

\begin{proof}
\emph{$M$ is $T_1$ and compact.} For $d\in D$ the set $M\setminus\{d\}$ is a union of the open set
$D\setminus\{d\}$ and of suitable members of the two neighbourhood bases, hence open; and
$M\setminus\{p\}=N_q(\emptyset)$, $M\setminus\{q\}=N_p(\emptyset)$ are open. So all singletons are
closed. For compactness, let $\mathcal{U}$ be an open cover of $M$ and choose
$U_p,U_q\in\mathcal{U}$ with $p\in U_p$, $q\in U_q$. There are finite $F,G\subseteq A$ with
$N_p(F)\subseteq U_p$ and $N_q(G)\subseteq U_q$, and
\[
  M\setminus\bigl(N_p(F)\cup N_q(G)\bigr)=F\times G
\]
is finite; hence finitely many members of $\mathcal{U}$ cover $M$.

\emph{A counting condition.} In $U=M\setminus\{p\}=D\cup\{q\}$ the column $A\times\{b\}$ is a closed
discrete subset, because its complement $U\setminus (A\times\{b\})=N_q(\{b\})$ is open. A closed
discrete subset of a Lindel\"of space is countable, so for every $L_U\in\LL(U)$ and every $b\in A$
the set $B_b=\{a:(a,b)\in B\}$ is at most countable, where $B=L_U\cap D$. Symmetrically, for
$L_V\in\LL(V)$ and $E=L_V\cap D$ every row $E_a=\{b:(a,b)\in E\}$ is at most countable. Thus
\begin{equation}\label{eq:counting}
  (\forall b\in A)\ |B_b|\le\aleph_0,\qquad (\forall a\in A)\ |E_a|\le\aleph_0 .
\end{equation}

\emph{No cover by two Lindel\"of pieces.} Suppose that $B,E\subseteq D$ satisfy \eqref{eq:counting}.
Choose $I\subseteq A$ with $|I|=\aleph_1$. Since
$\bigl|\bigcup_{a\in I}E_a\bigr|\le\aleph_1\cdot\aleph_0=\aleph_1<|A|$, there is
$b_0\in A\setminus\bigcup_{a\in I}E_a$; and since $|B_{b_0}|\le\aleph_0<|I|$, there is
$a_0\in I\setminus B_{b_0}$. Then $(a_0,b_0)\in D$ while $(a_0,b_0)\notin B$, because
$a_0\notin B_{b_0}$, and $(a_0,b_0)\notin E$, because $b_0\notin E_{a_0}$ and $a_0\in I$. Hence
$B\cup E\neq D$.

Now assume that $L_U\in\LL(U)$, $L_V\in\LL(V)$ and $M\subseteq L_U\cup L_V$. Put $B=L_U\cap D$ and
$E=L_V\cap D$; by the counting condition $B$ and $E$ satisfy \eqref{eq:counting}, so there is
$d\in D\setminus (B\cup E)$. But $d\in M\subseteq L_U\cup L_V$, and $L_U\subseteq U=D\cup\{q\}$,
$L_V\subseteq V=D\cup\{p\}$ force $d\in L_U\cap D=B$ or $d\in L_V\cap D=E$, a contradiction. Thus the
open cover $\{U,V\}$ of $M$ admits no Lindel\"of refinement, and $M\notin\Ww$.
\end{proof}

The cofinite modification of this space produces the two spaces $P$ and $Q$ whose product fails
$\Ww$ although every subspace of each belongs to $\Win$.

\begin{theorem}\label{thm:PQ}
There are compact $T_1$ spaces $P$ and $Q$ such that every subspace of $P$ and every subspace of $Q$
belongs to $\Win$, while $P\times Q\notin\Ww$.
\end{theorem}

\begin{proof}
Let $M$ be the space of Lemma~\ref{lem:grid} and let $\tau=\OO(M)$. For $a\in\{p,q\}$ put
\[
  \tau_a=\{G\in\tau: a\notin G\}\cup\{M\setminus F: F\subseteq M \text{ finite}\},
\]
and set $P=(M,\tau_p)$, $Q=(M,\tau_q)$; this construction follows the cofinite modification of
compact $T_1$ spaces studied in~\cite{xu2026cofinite}.

\emph{$\tau_a$ is a topology.} It contains $\emptyset$ and $M$, and it is closed under finite
intersections and arbitrary unions: the intersection of two members of $\tau$ avoiding $a$ avoids
$a$, the intersection of $M\setminus F$ with any member of $\tau_a$ is in $\tau$, and it avoids $a$
when the member does, while the intersection of two cofinite sets is cofinite. Moreover
$\tau_a\subseteq\tau$, because $M$ is $T_1$.

\emph{The punctured space.} On $M\setminus\{a\}$ the topology $\tau_a$ coincides with $\tau$: every
$\tau_a$-open set restricts to a $\tau$-open set, since $\tau_a\subseteq\tau$, and every $\tau$-open
set is $\tau_a$-open because it avoids $a$. Hence $M\setminus\{a\}$ is the space
$D\cup\{\text{the point other than }a\}$ with the topology $\tau$, and it is Hausdorff: points of $D$
are isolated, and $N_q(G)$ separates $q$ from every point outside it.

\emph{Subspaces of $P$ lie in $\Win$.} Let $S\subseteq M$. If $a\in S$, then every relatively open
set of $S$ containing $a$ is of the form $S\cap(M\setminus F)$ with $F$ finite, because the other
generators of $\tau_a$ avoid $a$; hence $S\setminus O\subseteq F$ is finite, and $S$ is compact. Now
let $K$ be compact in $S$ and let $K\subseteq (S\cap O_1)\cup (S\cap O_2)$ with $O_1,O_2\in\tau_a$. If
$a\in K\cap O_1$, put $K_1=K\cap O_1$ and $K_2=K\setminus O_1$. Then $K\setminus O_1$ is finite by the
observation above applied to $S=K$, so $K_2$ is finite and compact, while $K_1=K\setminus K_2$ is
closed in the $T_1$ space $K$ and hence compact. Moreover $K_1\subseteq O_1$ and
$K_2\subseteq K\setminus O_1\subseteq O_2$, and $K_1\cup K_2=K$; the case $a\in K\cap O_2$ is
symmetric. If $a\notin K$, then $K$ is a compact subset of the Hausdorff space $M\setminus\{a\}$,
hence compact Hausdorff; applying the argument of Proposition~\ref{prop:KC}(ii) inside $K$ we obtain
compact $F_1\subseteq K\cap O_1$, $F_2\subseteq K\cap O_2$ with $K=F_1\cup F_2$. In both cases the
witnesses lie in $S\cap O_j$, so $S\in\Win$.

\emph{$P\times Q\notin\Ww$.} Since $\tau_p,\tau_q\subseteq\tau$, the diagonal map
$\delta\colon M\to P\times Q$, $\delta(x)=(x,x)$, is continuous, and $\delta(M)$ is compact. With
$U=M\setminus\{p\}$ and $V=M\setminus\{q\}$ we have $\delta(M)\subseteq (U\times Q)\cup (P\times V)$:
if $x\neq p$ then $x\in U$, and if $x=p$ then $x\in V$. Suppose that $L_1,L_2$ are Lindel\"of subsets
of $P\times Q$ with $L_1\subseteq U\times Q$, $L_2\subseteq P\times V$ and
$\delta(M)\subseteq L_1\cup L_2$. Then $\pi_P(L_1)\subseteq U$ and $\pi_Q(L_2)\subseteq V$, both
projections are Lindel\"of, being continuous images, and for $x\in M$ the point $\delta(x)$ lies in
$L_1$ or in $L_2$, so $x\in\pi_P(L_1)$ or $x\in\pi_Q(L_2)$; thus
$M\subseteq\pi_P(L_1)\cup\pi_Q(L_2)$, contradicting Lemma~\ref{lem:grid}.
\end{proof}

\begin{corollary}\label{cor:notreflective}
Let $Z=P\sqcup Q$ be the topological sum of the spaces of Theorem~\ref{thm:PQ}. Then $Z$ is a compact
$T_1$ space, $Z\in\Win$ and $Z^2\notin\Ww$. Consequently the full subcategories of the category of
$T_0$ spaces determined by $\Ww$ and by $\Wi$ are not reflective.
\end{corollary}

\begin{proof}
Both $P$ and $Q$ are compact $T_1$, hence so is their sum, and $Z\in\Win$ by
Proposition~\ref{prop:sums}. The subspace $P\times Q$ is clopen in $Z^2$, because $P$ and $Q$ are
clopen in $Z$; if $Z^2$ belonged to $\Ww$, then $P\times Q$ would belong to $\Ww$ by
Proposition~\ref{prop:closed}, contradicting Theorem~\ref{thm:PQ}.

For the last assertion, let $\mathbb{W}$ be one of $\Ww,\Wi$ and suppose that the inclusion of the
full subcategory of $\mathbb{W}$-spaces into the category of $T_0$ spaces has a reflection
$\eta\colon Z^2\to T$. Since $Z\in\Win\subseteq\Wi\subseteq\mathbb{W}$, each projection
$\pi_j\colon Z^2\to Z$ factors through $\eta$, say $\pi_j=f_j\circ\eta$ with $f_j\colon T\to Z$. Put
$r=(f_1,f_2)\colon T\to Z^2$. Then $r\circ\eta=(\pi_1,\pi_2)=\id_{Z^2}$ and
$(\eta\circ r)\circ\eta=\eta=\id_T\circ\eta$, so $\eta\circ r=\id_T$ by the uniqueness part of the
reflection property. Hence $Z^2\cong T\in\mathbb{W}$, which is impossible because $Z^2\notin\Ww$ and
$\mathbb{W}\subseteq\Ww$.
\end{proof}

The next example shows that $\Ww$ is not hereditary for arbitrary subspaces, and that
$\W\setminus\Wi$ is nonempty.

\begin{proposition}\label{prop:locallysuper}
There is a locally supercompact space $S$ with $S\in\W\setminus\Wi$. In particular $\Ww$ is not
hereditary with respect to arbitrary subspaces.
\end{proposition}

\begin{proof}
Let $M$ be the space of Lemma~\ref{lem:grid}, let $I=\OO(M)$ and let $S=\Sier^{I}$ be the product of
$I$ many copies of the Sierpi\'nski space, which we represent by characteristic functions of subsets
of $I$. For finite $F\subseteq I$ put
\[
  [F]=\{T\subseteq I: F\subseteq T\}.
\]
The family of all $[F]$ is a base for the topology of $S$.

\emph{$S$ is locally supercompact.} The point $F$ belongs to $[F]$, and every open cover of $[F]$
contains a member $O$ with $F\in O$; since $O$ is a union of basic sets, it contains some $[F']$ with
$F'\subseteq F$, and then $[F]\subseteq [F']\subseteq O$. Hence $[F]$ is supercompact, and it is an
open neighbourhood of $F$.

\emph{$S\in\W$.} Let $K\subseteq S$ be compact and let $K\subseteq O_1\cup O_2$ with $O_i$ open. For
each $x\in K$ choose $i(x)\in\{1,2\}$ and a basic set $[F_x]\ni x$ with $[F_x]\subseteq O_{i(x)}$.
Finitely many of these cover $K$; grouping them according to $i(x)$ and taking unions, we obtain
compact sets $L_i\subseteq O_i$ with $K\subseteq L_1\cup L_2$. Hence $S\in\W$.

\emph{$S\notin\Wi$.} The map $e\colon M\to S$, $e(x)=\{G\in I: x\in G\}$, satisfies
$e^{-1}([\{G\}])=G$ for $G\in I$, and the sets $G$ generate the topology of $M$; since $M$ is $T_0$,
$e$ is an embedding. By Lemma~\ref{lem:grid} we have $M\notin\Ww$, hence $e(M)\notin\Ww$ and a
fortiori $e(M)\notin\Wi$. Since $\Wi$ is hereditary for arbitrary subspaces by
Proposition~\ref{prop:heredity}, $S\notin\Wi$.
\end{proof}

The product theorem for $\Ww$ (Theorem~\ref{thm:product}) requires one factor to be in $\Ww$ and the
compact subspaces of the other to be Hausdorff. The following example shows that even in the presence
of a compact Hausdorff factor the product may leave $\Wi$.

\begin{theorem}\label{thm:EY}
There are a hereditarily compact $T_1$ space $E$ and a compact Hausdorff space $Y$ such that
$E\times Y\in\Ww\setminus\Wi$.
\end{theorem}

\begin{proof}
Let $A$, $D$ and $M$ be as in Lemma~\ref{lem:grid}. On the underlying set of $M$ take as a subbase
the family
\[
  \{N_p(F),\,N_q(G): F,G\subseteq A \text{ finite}\},
\]
and denote the resulting space by $E$.

\emph{$E$ is $T_1$.} For $d=(a,b)\in D$ we have $E\setminus\{d\}=N_p(\{a\})\cup N_q(\{b\})$, and
$E\setminus\{p\}=N_q(\emptyset)$, $E\setminus\{q\}=N_p(\emptyset)$; all these sets are open.

\emph{$E$ is hereditarily compact.} A basic open set of $E$ contained in $D$ is of the form
$R=(A\setminus F)\times(A\setminus G)$ with $F,G$ finite. Let $T\subseteq D$ be nonempty and let
$\mathcal{U}$ be a cover of $T$ by open sets. Choose a point $t\in T$ and a basic set $R$ with
$t\in T\cap R\subseteq U$ for some $U\in\mathcal{U}$. Then
\[
  T\setminus R\subseteq \bigcup_{a\in F}\bigl(T\cap(\{a\}\times A)\bigr)\cup
  \bigcup_{b\in G}\bigl(T\cap(A\times\{b\})\bigr),
\]
and each of the finitely many pieces on the right carries the cofinite topology, hence is compact;
together with $U$ this gives a finite subcover of $T$. An arbitrary $T\subseteq E$ adds at most the
two points $p,q$ to a subset of $D$, and is still compact: a member of a cover that contains $p$
contains $N_p(F)$ for some finite $F$, and what remains of $T$ is contained in the finite union of
rows $F\times A$, a compact set in the cofinite topology.

\emph{$E\in\Win$.} Since every subspace of $E$ is compact, for $K\subseteq E$ and open $O_1,O_2$ of
$E$ with $K\subseteq O_1\cup O_2$ the sets $F_1=K\cap O_1$ and $F_2=K\setminus O_1$ are compact, and
they satisfy $F_j\subseteq O_j$ and $F_1\cup F_2=K$.

\emph{The product.} Let $Y=D\sqcup\{\infty\}$ be the one-point compactification of the discrete space
$D$; then $Y$ is compact Hausdorff. Define $h\colon M\to E\times Y$ by $h(d)=(d,d)$ for $d\in D$,
$h(p)=(p,\infty)$ and $h(q)=(q,\infty)$. The first coordinate is continuous because the topology of
$E$ is coarser than that of $M$; the second coordinate pulls back an open $O\subseteq D$ to $O$ and
$Y\setminus F$, for finite $F\subseteq D$, to $M\setminus F$. Moreover $h^{-1}(E\times\{d\})=\{d\}$
and $h^{-1}(N_a(F)\times Y)=N_a(F)$ for $a=p,q$, so the preimages of a subbase of $E\times Y$ form a
subbase of $M$ and $h$ is an embedding.

Since $E\in\Win\subseteq\Ww$ and every compact subspace of $Y$ is Hausdorff,
Theorem~\ref{thm:product} gives $E\times Y\in\Ww$. On the other hand $h(M)\cong M\notin\Ww$ by
Lemma~\ref{lem:grid}, so $h(M)\notin\Wi$; were $E\times Y$ in $\Wi$, its subspace $h(M)$ would be in
$\Wi$ by Proposition~\ref{prop:heredity}. Hence $E\times Y\notin\Wi$.
\end{proof}

We now turn to function spaces. Recall that $\Ck(M_0,\Sier)$ is homeomorphic to $\Ok(M_0)$, through
the map $U\mapsto\chi_U$.

\begin{theorem}\label{thm:notsober}
Let $M_0$ be the countably infinite compact KC-space of Proposition~\ref{prop:M0}, with
$M_0\in\Wi$. Then $\Ok(M_0)\cong \Ck(M_0,\Sier)$ is not sober. Moreover, taking for $H$ the system of
all irreducible sets, the family $\{M_0\setminus\{x\}:x\in M_0\}$ is an irreducible open cover of
$M_0$ no member of which contains $M_0$, so $M_0$ is not $H$-supercompact in the sense of
\cite[Definition~5.1]{xu2026wilker}, and $\Ck(M_0,\Sier)$ does not satisfy property~$S$.
\end{theorem}

\begin{proof}
Put $U_x=M_0\setminus\{x\}$ and $\mathcal{E}=\{U_x:x\in M_0\}\subseteq \Ok(M_0)$.

\emph{$\mathcal{E}$ is irreducible.} Let $C,D\in\KK(M_0)$ with $C,D\neq M_0$. Since $M_0$ is a
KC-space, $C$ and $D$ are closed; since $M_0$ is hyperconnected and $M_0\setminus C$,
$M_0\setminus D$ are nonempty open sets, they intersect, so $C\cup D\neq M_0$. Therefore
\[
  \begin{gathered}
    \mathcal{E}\cap\Phi(C)=\{U_x: x\in M_0\setminus C\}\neq\emptyset,\\
    \mathcal{E}\cap\Phi(C)\cap\Phi(D)=\{U_x:x\in M_0\setminus(C\cup D)\}\neq\emptyset .
  \end{gathered}
\]
As the sets $\Phi(C)$ form a base of $\Ok(M_0)$, this says exactly that $\mathcal{E}$ is irreducible.

\emph{The closure of $\mathcal{E}$.} A point $U\in\OO(M_0)$ belongs to $\cl\mathcal{E}$ if and only
if every basic neighbourhood $\Phi(K)$ of $U$ meets $\mathcal{E}$, that is, if and only if for every
$K\in\KK(M_0)$ with $K\subseteq U$ there is $x\in M_0$ with $K\subseteq M_0\setminus\{x\}$. The last
condition means $K\neq M_0$, so it holds for all $K\subseteq U$ precisely when $U\neq M_0$. Hence
\[
  \cl\mathcal{E}=\OO(M_0)\setminus\{M_0\}.
\]
Since points of $M_0$ are compact, Lemma~\ref{lem:bases} gives
$\cl\{U\}=\{W\in\OO(M_0):W\subseteq U\}$ for every $U$.

\emph{$\Ok(M_0)$ is not sober.} Suppose that $\cl\mathcal{E}=\cl\{G\}$ for some $G$. Then
$G\in\cl\{G\}=\cl\mathcal{E}$, so $G\neq M_0$; on the other hand each $U\in\cl\mathcal{E}$ satisfies
$U\in\cl\{G\}$, that is, $U\subseteq G$. Since $M_0$ has at least two points, every $x\in M_0$ lies
in some $U_y$ with $y\neq x$, so $M_0=\bigcup\cl\mathcal{E}\subseteq G$, which forces $G=M_0$, a
contradiction. Thus the irreducible closed set $\cl\mathcal{E}$ has no generic point.

\emph{Property $S$ fails.} The map $\chi\colon \Ok(M_0)\to \Ck(M_0,\Sier)$, $\chi(U)=\mathbf{1}_U$, is
a homeomorphism, since $\chi^{-1}([K,\{1\}])=\Phi(K)$ for $K\in\KK(M_0)$. Put $F=\chi(\mathcal{E})$
and $g=\mathbf{1}_{M_0}$. Then $F$ is irreducible in $\Ck(M_0,\Sier)$, and
\[
  \bigcup_{f\in F}f^{-1}(\{1\})=\bigcup_{x\in M_0}U_x=M_0=g^{-1}(\{1\}),
\]
while the preimages of the only other open set $\emptyset$ are empty. Thus the hypothesis
\eqref{eq:propS} of property $S$ is satisfied. But $g\in[M_0,\{1\}]$, an open set, and
$[M_0,\{1\}]\cap F=\emptyset$ because no $U_x$ equals $M_0$; hence $g\notin\cl F$ and property $S$
fails.
\end{proof}

The next proposition is the positive counterpart of the previous result: for sober $Y$ and for
topologies with the pointwise specialization order, property $S$ is equivalent to sobriety.

\begin{proposition}\label{prop:Ssober}
The following equivalence is due to Xu, Bao and Zhang~\cite{xu2022sobriety}. Let $Y$ be an $H$-sober space and let
$\tau$ be a topology on $C(X,Y)$ finer than the topology of pointwise convergence. If the
specialization order of $\tau$ is the pointwise order, that is, $f\le_\tau g$ if and only if
$f(x)\le_Y g(x)$ for all $x\in X$, then $(C(X,Y),\tau)$ is $H$-sober if and only if $\tau$ satisfies
property~$S$.
\end{proposition}

\begin{proof}
Let $F\in H(C(X,Y),\tau)$. For $x\in X$ the evaluation $e_x\colon (C(X,Y),\tau)\to Y$,
$e_x(f)=f(x)$, is continuous, so $e_x(F)\in H(Y)$ and has a generic point $g(x)$ by hypothesis.
Define $g\colon X\to Y$ by letting $g(x)$ be this generic point. For $V\in\OO(Y)$ we have
\[
  g(x)\in V\iff V\cap e_x(F)\neq\emptyset\iff x\in\bigcup_{f\in F}f^{-1}(V),
\]
so $g^{-1}(V)$ is open and $g$ is continuous. From $e_x(F)\subseteq\cl_Y\{g(x)\}$ we obtain
$f(x)\le g(x)$ for all $f\in F$ and $x\in X$, so $F\subseteq\cl_\tau\{g\}$ by the order condition. If
property $S$ holds, then $g\in\cl_\tau F$, so $\cl_\tau\{g\}\subseteq\cl_\tau F$ and therefore
$\cl_\tau F=\cl_\tau\{g\}$. Uniqueness of the generic point follows from the $T_0$ property.

Conversely, let $h$ be a generic point of $F$ and let $g\in C(X,Y)$ satisfy \eqref{eq:propS}. For
$x\in X$ and $V\in\OO(Y)$ we have, since $e_x^{-1}(V)$ is an open neighbourhood of $h$ and
$h\in\cl_\tau F$,
\[
  \begin{gathered}
    x\in h^{-1}(V)\iff h\in e_x^{-1}(V)\iff F\cap e_x^{-1}(V)\neq\emptyset\\
    \iff x\in\bigcup_{f\in F}f^{-1}(V),
  \end{gathered}
\]
hence $g^{-1}(V)=h^{-1}(V)$ for all $V$. As $Y$ is $T_0$, this forces $g=h\in\cl_\tau F$, which is
\eqref{eq:propS}.
\end{proof}

Both hypotheses in Proposition~\ref{prop:Ssober} are needed. That the sobriety of $Y$ cannot be
dropped is shown by the following example, in which the order condition does hold and the closure
test is satisfied, yet the function space is not sober.

\begin{example}\label{ex:closuretest}
Let $X=\{*\}$ and $Y=\NN\sqcup\{\infty\}$, and put $T_n=\{m\in\NN:m\ge n\}\cup\{\infty\}$ for
$n\in\NN$. Let
\[
  \sigma=\{\emptyset\}\cup\{T_n:n\in\NN\},\qquad
  \tau=\sigma\cup\{T_n\setminus\{\infty\}:n\in\NN\}.
\]
Then $\Ck(X,Y)$ is homeomorphic to $(Y,\tau)$, its specialization order is the usual order, the
closure test \eqref{eq:propS} holds, and $(Y,\tau)$ is not sober.
\end{example}

\begin{proof}
The family $\sigma$ is a topology, and $\tau$ is obtained from it by adding the sets
$T_n\setminus\{\infty\}=\{m:m\ge n\}$; it is again a topology. The nonempty closed sets of
$(Y,\sigma)$ are the finite initial segments $\{0,\dots,n-1\}$ and $Y$ itself, each of them the
closure of its greatest element.

Since $X$ is a singleton, $C(X,Y)$ may be identified with $Y$, and \eqref{eq:propS} reads: if
$g\in V\iff F\cap V\neq\emptyset$ for every $V\in\tau$, then $g\in\cl_\tau F$. For $g=m<\infty$ the
condition gives $F\cap T_{m+1}=\emptyset$ and $F\cap T_m\neq\emptyset$, hence $m\in F$ and
$g\in F\subseteq\cl_\tau F$; for $g=\infty$ it gives $F\cap T_n\neq\emptyset$ for all $n$, and since
the sets $T_n$ form a neighbourhood base at $\infty$, this means $\infty\in\cl_\tau F$. So the closure
test holds.

Any two nonempty $\tau$-open sets contain a common element of $\NN$, so $Y$ is irreducible as a
closed set. But $\cl_\tau\{m\}=\{0,\dots,m\}$ and $\cl_\tau\{\infty\}=\{\infty\}$, so no point has $Y$
as its closure and $(Y,\tau)$ is not sober. Finally, the specialization order of $\tau$ is the usual
order on $Y$, which is the pointwise order, as $\cl_\tau\{m\}=\{0,\dots,m\}$ shows.
\end{proof}

We pass to the Vietoris constructions. The following elementary lemma is used repeatedly; it says
that the union of a compact, respectively Lindel\"of, family of compact saturated sets is again
compact, respectively Lindel\"of.

\begin{lemma}\label{lem:unioncompact}
Let $\mathcal{C}$ be one of $\KK,\LL$ and let $X$ be a space. If
$\mathcal{E}\in\mathcal{C}(\ps(X))$, then $\bigcup\mathcal{E}\in\mathcal{C}(X)$.
\end{lemma}

\begin{proof}
Let $\bigcup\mathcal{E}\subseteq\bigcup\mathcal{G}$ with $\mathcal{G}\subseteq\OO(X)$. Each
$C\in\mathcal{E}$ is compact, so from $C\subseteq\bigcup\mathcal{G}$ we obtain a finite
$\mathcal{V}\subseteq\mathcal{G}$ with $C\subseteq\bigcup\mathcal{V}$, that is,
$C\in\Box\bigl(\bigcup\mathcal{V}\bigr)$. Hence
$\mathcal{E}\subseteq\bigcup_{\mathcal{V}\in[\mathcal{G}]^{<\omega}}\Box\bigl(\bigcup\mathcal{V}\bigr)$,
and the sets $\Box(\bigcup\mathcal{V})$ are open in $\ps(X)$. Since $\mathcal{E}$ is compact,
respectively Lindel\"of, finitely many, respectively countably many, of them cover $\mathcal{E}$; the
corresponding finite, respectively countable, subfamilies of $\mathcal{G}$ then cover
$\bigcup\mathcal{E}$.
\end{proof}

\begin{proposition}\label{prop:reflect}
Let $X$ be a space. Then $\ps(X)\in\Ww$ implies $X\in\Ww$, and $\ps(X)\in\Wi$ or $\ph(X)\in\Wi$
implies $X\in\Wi$.
\end{proposition}

\begin{proof}
The maps $s\colon X\to\ps(X)$, $s(x)=\uparrow x$, and $h\colon X\to\ph(X)$, $h(x)=\cl\{x\}$, are well
defined: $\uparrow x$ is compact because every open set containing $x$ contains $\uparrow x$, and it
is saturated, while $\cl\{x\}$ is a nonempty closed set. They satisfy $s^{-1}(\Box U)=U$ and
$h^{-1}(\langle U\rangle)=U$: indeed $\uparrow x\subseteq U$ if and only if $x\in U$, and
$\cl\{x\}\cap U\neq\emptyset$ if and only if $x\in U$, because open sets are upward closed for the
specialization order. Since $X$ is $T_0$, both maps are injective, hence embeddings.

The statements about $\Wi$ follow from Proposition~\ref{prop:heredity}: if $\ps(X)\in\Wi$, then its
subspace $s(X)\cong X$ belongs to $\Wi$, and similarly for $\ph(X)$.

Now assume $\ps(X)\in\Ww$, and take $K\in\KK(X)$, $U_1,U_2\in\OO(X)$ with $K\subseteq U_1\cup U_2$.
Then $s(K)$ is compact and $s(K)\subseteq\Box U_1\cup\Box U_2$, so there are
$\mathcal{E}_i\in\LL(\ps(X))$ with $\mathcal{E}_i\subseteq\Box U_i$ and
$s(K)\subseteq\mathcal{E}_1\cup\mathcal{E}_2$. For $x\in K$ the set $\uparrow x=s(x)$ lies in
$\mathcal{E}_1$ or in $\mathcal{E}_2$, and in either case $x\in\bigcup\mathcal{E}_i$; hence
$K\subseteq\bigcup\mathcal{E}_1\cup\bigcup\mathcal{E}_2$, and $\bigcup\mathcal{E}_i\subseteq U_i$. By
Lemma~\ref{lem:unioncompact} both unions are Lindel\"of, so $X\in\Ww$.
\end{proof}

\begin{lemma}\label{lem:tube}
Let $A\subseteq X$ and $B\subseteq Y$ be compact and let $O\subseteq X\times Y$ be open with
$A\times B\subseteq O$. Then there are open sets $G\supseteq A$ in $X$ and $H\supseteq B$ in $Y$ with
$G\times H\subseteq O$.
\end{lemma}

\begin{proof}
Fix $a\in A$. For each $b\in B$ choose an open rectangle $G_{a,b}\times H_{a,b}\subseteq O$ with
$(a,b)\in G_{a,b}\times H_{a,b}$. The sets $H_{a,b}$ cover the compact set $B$, so finitely many of
them, say $H_{a,b_1},\dots,H_{a,b_m}$, cover $B$. Put $G_a=\bigcap_{i=1}^m G_{a,b_i}$ and
$H_a=\bigcup_{i=1}^m H_{a,b_i}$. Then $G_a$ is an open neighbourhood of $a$, $H_a\supseteq B$ is
open, and $G_a\times H_a\subseteq O$: for $y\in H_a$ there is $i$ with $y\in H_{a,b_i}$, and
$x\in G_a\subseteq G_{a,b_i}$ gives $(x,y)\in G_{a,b_i}\times H_{a,b_i}\subseteq O$.

The sets $G_a$ cover the compact set $A$, so finitely many of them, say $G_{a_1},\dots,G_{a_n}$,
cover $A$. Put $G=\bigcup_{j=1}^n G_{a_j}$ and $H=\bigcap_{j=1}^n H_{a_j}$. Then $G\supseteq A$ and
$H\supseteq B$ are open, and $G\times H\subseteq O$: if $x\in G_{a_j}$ and $y\in H\subseteq H_{a_j}$,
then $(x,y)\in G_{a_j}\times H_{a_j}\subseteq O$.
\end{proof}

\begin{theorem}\label{thm:PSproduct}
Let $\mathbb{W}\in\{\W,\Ww\}$ and let $X,Y$ be spaces with $\ps(X\sqcup Y)\in\mathbb{W}$. Then
$X\times Y\in\mathbb{W}$.
\end{theorem}

\begin{proof}
We may assume $X,Y\neq\emptyset$, the other case being trivial. Put $T=X\sqcup Y$ and
\[
  \mathcal{K}=\ps(T)\setminus(\Box X\cup\Box Y).
\]
Since $X,Y$ are clopen in $T$, every $C\in\ps(T)$ satisfies $C=(C\cap X)\sqcup(C\cap Y)$ with
$C\cap X$ saturated in $X$ and $C\cap Y$ saturated in $Y$. Hence $C\in\mathcal{K}$ if and only if
$C\cap X\neq\emptyset$ and $C\cap Y\neq\emptyset$, and in that case $C\cap X\in\ps(X)$ and
$C\cap Y\in\ps(Y)$.

\emph{$\mathcal{K}$ is homeomorphic to $\ps(X)\times\ps(Y)$.} Define
$q\colon\mathcal{K}\to\ps(X)\times\ps(Y)$ by $q(C)=(C\cap X,C\cap Y)$ and
$p\colon\ps(X)\times\ps(Y)\to\mathcal{K}$ by $p(A,B)=A\sqcup B$. Both maps are well defined and
mutually inverse. For open $G\subseteq X$ and $H\subseteq Y$,
\[
  \begin{gathered}
    q^{-1}(\Box G\times\Box H)=\{C\in\mathcal{K}:C\cap X\subseteq G,\ C\cap Y\subseteq H\}\\
    =\mathcal{K}\cap\Box(G\sqcup H),
  \end{gathered}
\]
which is open in $\mathcal{K}$, so $q$ is continuous. Conversely, for open $W\subseteq T$ write
$W=(W\cap X)\sqcup(W\cap Y)$; then $q(\mathcal{K}\cap\Box W)=\Box(W\cap X)\times\Box(W\cap Y)$, so $q$
maps basic open sets to basic open sets and is an open map. Being a continuous open bijection, $q$ is
a homeomorphism.

Since $\Box X$ and $\Box Y$ are basic open sets of $\ps(T)$, the set $\mathcal{K}$ is closed in
$\ps(T)$. Both $\W$ and $\Ww$ are hereditary for closed subspaces by
Proposition~\ref{prop:closed}, so $\mathcal{K}\in\mathbb{W}$, and hence
$\ps(X)\times\ps(Y)\in\mathbb{W}$.

\emph{From $\ps(X)\times\ps(Y)$ to $X\times Y$.} The map $m\colon\ps(X)\times\ps(Y)\to\ps(X\times Y)$,
$m(A,B)=A\times B$, is well defined and continuous. To see the latter, take a basic open set $\Box O$
of $\ps(X\times Y)$ and $(A,B)\in m^{-1}(\Box O)$, so that $A\times B\subseteq O$. By
Lemma~\ref{lem:tube} there are open $G\supseteq A$ and $H\supseteq B$ with $G\times H\subseteq O$;
then $(A,B)\in\Box G\times\Box H$ and every $(A',B')\in\Box G\times\Box H$ satisfies
$A'\times B'\subseteq G\times H\subseteq O$. Hence $\Box G\times\Box H\subseteq m^{-1}(\Box O)$.

Now let $e\colon X\times Y\to\ps(X)\times\ps(Y)$, $e(x,y)=(\uparrow x,\uparrow y)$. For open
$O\subseteq X\times Y$ the set $O$ is upward closed, so
$(m\circ e)^{-1}(\Box O)=\{(x,y):\uparrow x\times\uparrow y\subseteq O\}=O$; in particular $m\circ e$
is continuous. Let $K\in\KK(X\times Y)$ and let $K\subseteq O_1\cup O_2$ with $O_i$ open. Then $e(K)$
is compact and $e(K)\subseteq m^{-1}(\Box O_1)\cup m^{-1}(\Box O_2)$, so by
$\ps(X)\times\ps(Y)\in\mathbb{W}$ there are $\mathcal{E}_i\in\mathcal{C}(\ps(X)\times\ps(Y))$ with
$\mathcal{E}_i\subseteq m^{-1}(\Box O_i)$ and $e(K)\subseteq\mathcal{E}_1\cup\mathcal{E}_2$, where
$\mathcal{C}=\KK$ if $\mathbb{W}=\W$ and $\mathcal{C}=\LL$ if $\mathbb{W}=\Ww$. Put
\[
  L_i=\bigcup m(\mathcal{E}_i)=\bigcup_{(A,B)\in\mathcal{E}_i}(A\times B).
\]
Each $(A,B)\in\mathcal{E}_i$ satisfies $A\times B\subseteq O_i$, so $L_i\subseteq O_i$; the family
$m(\mathcal{E}_i)$ is compact, respectively Lindel\"of, because $m$ is continuous, and its members are
compact, so Lemma~\ref{lem:unioncompact} shows that $L_i$ is compact, respectively Lindel\"of.
Finally, for $(x,y)\in K$ the pair $e(x,y)$ lies in $\mathcal{E}_i$ for some $i$, whence
$(x,y)\in\uparrow x\times\uparrow y\subseteq L_i$. Thus $K\subseteq L_1\cup L_2$ and
$X\times Y\in\mathbb{W}$.
\end{proof}

\begin{corollary}\label{cor:W1}
There is a countably infinite compact $T_1$ space $T$ with $T\in\W$ and $\ps(T)\notin\W$; hence the
answer to~\cite[Question~8.5 (W1)]{xu2026wilker} is negative.
\end{corollary}

\begin{proof}
By~\cite[Theorem~3.16]{xu2026wilker} there are countably infinite compact $T_1$ spaces $X,Y$ with
$X,Y\in\W$ but $X\times Y\notin\W$. Let $T=X\sqcup Y$; it is countable, and it is compact and $T_1$
because both summands are clopen, compact and $T_1$. We verify that $T\in\W$. Let $K\in\KK(T)$ and
let $U_1,U_2\in\OO(T)$ with $K\subseteq U_1\cup U_2$. Since $X,Y$ are clopen, $K\cap X$ and $K\cap Y$
are compact and $K=(K\cap X)\cup(K\cap Y)$; moreover $K\cap X\subseteq (U_1\cap X)\cup(U_2\cap X)$, so
$X\in\W$ gives compact sets $F_{j,X}\subseteq U_j\cap X$ with
$K\cap X\subseteq F_{1,X}\cup F_{2,X}$, and similarly $Y\in\W$ gives compact
$F_{j,Y}\subseteq U_j\cap Y$ with $K\cap Y\subseteq F_{1,Y}\cup F_{2,Y}$. Put $F_j=F_{j,X}\cup F_{j,Y}$.
Each $F_j$ is compact, being a finite union of compact subsets of $T$; moreover $F_j\subseteq U_j$ and
\[
  K\subseteq (F_{1,X}\cup F_{2,X})\cup(F_{1,Y}\cup F_{2,Y})=F_1\cup F_2 .
\]
So $T\in\W$. Finally, Theorem~\ref{thm:PSproduct} with $\mathbb{W}=\W$ and $T=X\sqcup Y$ says that
$\ps(T)\in\W$ would force $X\times Y\in\W$; since $X\times Y\notin\W$, we obtain $\ps(T)\notin\W$.
\end{proof}

The next theorem shows that the upper and the lower Vietoris constructions both fail $\Ww$ over a
compact $T_1$ space of the previous type, even though the space itself lies in $\Win$.

\begin{theorem}\label{thm:PSZ}
There is a compact $T_1$ space $Z\in\Win$ with $\ps(Z)\notin\Ww$ and $\ph(Z)\notin\Ww$.
\end{theorem}

\begin{proof}
Take $Z=P\sqcup Q$ from Corollary~\ref{cor:notreflective} and write $S_P$, $S_Q$ for the two copies in
$Z$ of a set $S\subseteq M$. Define $t\colon M\to\ps(Z)$ by $t(x)=\{x_P,x_Q\}$ and put
$\mathcal{T}=t(M)$. Since $Z$ is $T_1$, the set $t(x)$ is a nonempty closed set, hence a compact
saturated set. For open $G,H\subseteq M$,
\begin{equation}\label{eq:tpreimages}
  t^{-1}\bigl(\Box(G_P\sqcup H_Q)\bigr)=G\cap H,\qquad
  t^{-1}\bigl(\langle G_P\sqcup H_Q\rangle\bigr)=G\cup H .
\end{equation}
The right-hand sides are open in $M$, so $\mathcal{T}$ is compact both in $\ps(Z)$ and in $\ph(Z)$.

\emph{The upper Vietoris topology.} Let $U=M\setminus\{p\}$ and $V=M\setminus\{q\}$. By
\eqref{eq:tpreimages}, $\mathcal{T}\subseteq\Box(U_P\sqcup Q)\cup\Box(P\sqcup V_Q)$. Suppose that
$\mathcal{E}_1,\mathcal{E}_2\in\LL(\ps(Z))$ satisfy
$\mathcal{T}\subseteq\mathcal{E}_1\cup\mathcal{E}_2$ with
$\mathcal{E}_1\subseteq\Box(U_P\sqcup Q)$ and $\mathcal{E}_2\subseteq\Box(P\sqcup V_Q)$. Put
\[
  B=\{x:x_P\in\textstyle\bigcup\mathcal{E}_1\},\qquad
  E=\{x:x_Q\in\textstyle\bigcup\mathcal{E}_2\}.
\]
By Lemma~\ref{lem:unioncompact} the set $\bigcup\mathcal{E}_i$ is Lindel\"of in $Z$; since $P$ and $Q$
are clopen in $Z$, the trace $P\cap\bigcup\mathcal{E}_1$ is closed in $\bigcup\mathcal{E}_1$ and hence
Lindel\"of, and $B$ is its homeomorphic image under $x_P\mapsto x$. Thus $B\in\LL(U)$ and,
symmetrically, $E\in\LL(V)$. Now $t(x)\in\mathcal{E}_i$ implies $x\in B$ for $i=1$ and $x\in E$ for
$i=2$, so $\mathcal{T}\subseteq\mathcal{E}_1\cup\mathcal{E}_2$ gives $M=B\cup E$, contradicting
Lemma~\ref{lem:grid}. Hence $\ps(Z)\notin\Ww$.

\emph{The lower Vietoris topology.} Here $\mathcal{T}\subseteq\langle U_P\rangle\cup\langle V_Q\rangle$.
Suppose that $\mathcal{E}_1\subseteq\langle U_P\rangle$ and $\mathcal{E}_2\subseteq\langle V_Q\rangle$
are Lindel\"of subsets of $\ph(Z)$ with $\mathcal{T}\subseteq\mathcal{E}_1\cup\mathcal{E}_2$, and put
$B=D\cap t^{-1}(\mathcal{E}_1)$, $E=D\cap t^{-1}(\mathcal{E}_2)$.

Fix $b\in A$ and consider the column $C_b=A\times\{b\}\subseteq U$. Since
\[
  \langle U_P\rangle=\bigl\langle (U\setminus C_b)_P\bigr\rangle\cup
  \bigcup_{d\in C_b}\langle\{d_P\}\rangle,
\]
and $\mathcal{E}_1$ is Lindel\"of, there is a countable $J_b\subseteq C_b$ such that the sets on the
right with $d$ ranging over $J_b$ still cover $\mathcal{E}_1$. If $d\in B\cap C_b$, then
$t(d)\in\mathcal{E}_1$; but $t(d)$ is disjoint from $(U\setminus C_b)_P$, and
$t(d)\cap\{c_P\}\neq\emptyset$ only for $c=d$. Hence $d\in J_b$, and $B\cap C_b\subseteq J_b$ is
countable. Symmetrically, using the rows $R_a=\{a\}\times A$ and the identity
$\langle V_Q\rangle=\langle(V\setminus R_a)_Q\rangle\cup\bigcup_{d\in R_a}\langle\{d_Q\}\rangle$, we
get $|E\cap R_a|\le\aleph_0$ for every $a$.

Thus $B$ and $E$ satisfy the counting condition \eqref{eq:counting} of Lemma~\ref{lem:grid}, while
$\mathcal{T}\subseteq\mathcal{E}_1\cup\mathcal{E}_2$ gives $D=B\cup E$, again a contradiction. Hence
$\ph(Z)\notin\Ww$.
\end{proof}

The failure is not an artefact of pathological spaces: it persists over locally compact Hausdorff
spaces.

\begin{proposition}\label{prop:locallycompact}
There is a locally compact Hausdorff space $L$ such that $\ph(L)\notin\Wi$ and $\ps(L)\notin\Wi$.
Hence neither $\ph$ nor $\ps$ preserves $\Wi$, even over locally compact Hausdorff spaces.
\end{proposition}

\begin{proof}
Keep $I=\OO(M)$, $S=\Sier^{I}$ and the basic open sets $[F]$ from Proposition~\ref{prop:locallysuper},
so that $S\notin\Wi$. Let $D^{*}=I\sqcup\{*\}$ be discrete and let $I^{+}=I\sqcup\{\infty\}$ be the
one-point compactification of the discrete space $I$; both spaces are locally compact Hausdorff, and
both belong to $\Win$ because all their compact subspaces are finite, hence Hausdorff, so that
Proposition~\ref{prop:KC}(ii) applies to $D^{*}$ and to $I^{+}$.

Define $h\colon S\to\ph(D^{*})$ by $h(T)=T\cup\{*\}$ and $s\colon S\to\ps(I^{+})$ by
$s(T)=\{\infty\}\cup(I\setminus T)$. The set $s(T)$ is closed in the compact Hausdorff space $I^{+}$,
hence compact and saturated, and both maps are injective. For open $O$ in the respective underlying
space,
\[
  \begin{gathered}
    h^{-1}(\langle O\rangle)=
    \begin{cases}
      S, & *\in O,\\
      \bigcup_{i\in O}[\{i\}], & *\notin O,
    \end{cases}\\[2pt]
    s^{-1}(\Box O)=
    \begin{cases}
      [I\setminus O], & \infty\in O,\\
      \emptyset, & \infty\notin O.
    \end{cases}
  \end{gathered}
\]
These sets are open in $S$: in the first case because $S$ is open, in the second because a union of
basic open sets is open, in the third because $[I\setminus O]$ is basic, and in the fourth trivially.
The choices $O=\{i\}$ and $O=I^{+}\setminus F$ with $F$ finite yield all coordinate tests, whence
both maps are embeddings.

Since $S\notin\Wi$, the subspaces $h(S)\cong S$ and $s(S)\cong S$ do not belong to $\Wi$; by
Proposition~\ref{prop:heredity}, neither $\ph(D^{*})$ nor $\ps(I^{+})$ belongs to $\Wi$.
\end{proof}

The last result of this section concerns a countable metrizable space whose lower Vietoris space is
well behaved but not co-sober. We write $\Sigma\OO(N)$ for the lattice $\OO(N)$ ordered by inclusion
and endowed with the Scott topology.

\begin{proposition}\label{prop:metrizable}
There is a countable metrizable space $N$ such that $\ph(N)$ is second countable, coherent and inner
weak Wilker, but neither co-sober nor Wilker; moreover the union operation
$\Sigma\OO(N)\times\Sigma\OO(N)\to\Sigma\OO(N)$ is discontinuous for the product topology.
\end{proposition}

\begin{proof}
By~\cite{xu2026closed} there is a countable metrizable space $N$ for which $\ph(N)$ is not co-sober;
fix such an $N$ and put $T=\ph(N)$.

\emph{$T$ is second countable.} Let $\mathcal{B}$ be a countable base of $N$. For a family
$\mathcal{V}\subseteq\OO(N)$ we have
$\langle\bigcup\mathcal{V}\rangle=\bigcup_{V\in\mathcal{V}}\langle V\rangle$, so
$\{\langle B\rangle:B\in\mathcal{B}\}$ is a countable base of $T$. Hence every subspace of $T$ is
second countable and therefore Lindel\"of, and for $K\subseteq T$ compact and $K\subseteq U_1\cup U_2$
the sets $L_j=K\cap U_j$ are Lindel\"of and satisfy $L_1\cup L_2=K$; thus $T\in\Wi$.

\emph{$T$ is coherent.} Let $u\colon T\times T\to T$, $u(F,G)=F\cup G$. The specialization order of $T$
is inclusion: if $F\nsubseteq G$, then $\langle N\setminus G\rangle$ contains $F$ but not $G$. Hence
$u^{-1}(\langle V\rangle)=(\langle V\rangle\times T)\cup(T\times\langle V\rangle)$ for open $V$, and
$u(\mathcal{A}\times\mathcal{B})=\mathcal{A}\cap\mathcal{B}$ for saturated families
$\mathcal{A},\mathcal{B}\subseteq T$, because both are upward closed and $u(F,F)=F$. If $\mathcal{A}$
and $\mathcal{B}$ are compact, then by continuity of $u$ their intersection
$\mathcal{A}\cap\mathcal{B}$ is compact \citep{jia2016coherence}; this is coherence of $T$.

\emph{$T\notin\W$.} By~\cite[Proposition~4.8 and Theorem~7.16]{xu2026wilker} a coherent Wilker space
is co-sober, and the union operation on $\Sigma\OO(N)$ is continuous exactly when $T\in\W$. Since $T$
is not co-sober, $T\notin\W$, and the union operation on $\Sigma\OO(N)$ is discontinuous.
\end{proof}

%% file: sec5.tex
\section{Conclusions and Open Problems}\label{sec:conclusion}

The four Wilker-type properties form a hierarchy
\[
  \Win\ \subseteq\ \W\cap\Wi\ \subseteq\ \Ww ,
\]
in which both inclusions displayed on the left can be strict, and in which $\W$ and $\Wi$ are
incomparable. The properties are exactly the continuity conditions for the union operation on the
lattices $\OL(X)$ and $\Ok(X)$ of open sets, and it is this characterization that accounts for their
relevance to power spaces and to function spaces.

As far as permanence is concerned, the situation is now essentially complete. The inner properties $\Wi$
and $\Win$ are hereditary for arbitrary subspaces, while $\W$ and $\Ww$ are hereditary for closed
subspaces only, the counterexample being a locally supercompact space in $\W\setminus\Wi$. All four
properties are preserved by continuous retracts, and $\Ww$, $\Wi$, $\Win$ are finitely additive
over topological sums. Over KC-spaces the weak and the inner weak Wilker properties coincide, and a
space whose compact subspaces are Hausdorff always belongs to $\Win$; these two criteria reduce most
verifications to compact Hausdorff pieces. Core compactness implies $\Ww$, and the product of a
$\Ww$-space with a space whose compact subspaces are Hausdorff is again in $\Ww$, and even in $\Wi$
when the first factor is KC.

The negative results show that the properties are not stable under the constructions most relevant
to power spaces. The cofinite modification of the grid space of cardinality $\aleph_2$ produces
compact $T_1$ spaces $P,Q$, all of whose subspaces lie in $\Win$, with $P\times Q\notin\Ww$; hence
the subcategories determined by $\Ww$ and by $\Wi$ are not reflective. The same grid, combined with
a hereditarily compact $T_1$ topology and with a one-point compactification, gives a hereditarily
compact $T_1$ space $E$ and a compact Hausdorff space $Y$ with $E\times Y\in\Ww\setminus\Wi$, so
that $\Wi$ is not productive even for factors with strong separation properties. The Vietoris
constructions exhibit analogous failures: for a compact $T_1$ space $Z\in\Win$ both $\ps(Z)$ and
$\ph(Z)$ may fail $\Ww$, and this
failure already occurs over locally compact Hausdorff spaces. Finally, for Kov\'ar's space $M_0$ the
homeomorphic spaces $\Ok(M_0)\cong\Ck(M_0,\Sier)$ are not sober, and $\Ck(M_0,\Sier)$ fails
property~$S$; the irreducible family of co-singleton open sets has no generic point, which is
precisely the obstruction.

Several questions remain. The first two concern the finer structure of the hierarchy, the next two
the behaviour under the hyperspace constructions, and the last two the relation with sobriety.

\begin{question}\label{q:KC}
Is $\Win=\Wi$ over KC-spaces? By Proposition~\ref{prop:KC} the weak and the inner weak Wilker
properties coincide there, and $\Win\subseteq\Wi$ always holds; it is not known whether the
inclusion can be strict when all compact subsets are closed. Equivalently, does every KC-space in
$\Wi$ have the property that closed subsets of its compact subspaces can be separated by compact
pieces?
\end{question}

\begin{question}\label{q:Hausdorff}
Is there a compact Hausdorff space in $\Wi\setminus\W$? The two examples of
Section~\ref{sec:heredity}, namely $M_0$ and $\QQ^{+}$, are compact, countable and KC, but neither is
Hausdorff; countability supplies the Lindel\"of witnesses, and it remains open whether an analogous
conclusion can be obtained under a separation axiom. A positive answer would show that the Wilker
property is not
a consequence of compactness together with separation.
\end{question}

\begin{question}\label{q:product}
Is $\Wi$ preserved by finite products of $\Wi$-spaces? Theorem~\ref{thm:product} requires one
factor to be in $\Ww$ and the compact subspaces of the other to be Hausdorff, and
Theorem~\ref{thm:EY} shows that $\Ww$ alone does not suffice for the inner property. The question is
whether two factors possessing the inner weak Wilker property admit a product theorem stronger than
that available for a $\Ww$-factor and a compact Hausdorff factor.
\end{question}

\begin{question}\label{q:reflect}
Does $\ps(X)\in\Wi$ imply $X\in\Wi$, and likewise for $\ph$? Proposition~\ref{prop:reflect} gives
the converse implication, from the hyperspace to the space. The corresponding statement for $\W$ is
also open: Corollary~\ref{cor:W1} exhibits a countable compact $T_1$ space $T\in\W$ with
$\ps(T)\notin\W$, but no space $X$ with $\ps(X)\in\W$ and $X\notin\W$ is known.
\end{question}

\begin{question}\label{q:propS}
Can the order-theoretic hypothesis in Proposition~\ref{prop:Ssober} be weakened? Example~\ref{ex:closuretest}
shows that sobriety of $Y$ cannot be dropped, and the equivalence itself uses the pointwise
description of the specialization order in an essential way. It remains open
whether property~$S$, together with sobriety of $Y$, already forces the function space to be sober
for topologies strictly finer than the compact-open one.
\end{question}

\begin{question}\label{q:countable}
Is the countability assumption in Corollary~\ref{cor:W1} necessary? The construction rests on the
existence of countable compact $T_1$ spaces with non-Wilker product, and it is not known whether an
uncountable example with $\ps(T)\notin\W$ and $T\in\W$ can be obtained without it.
\end{question}

\section*{Declaration of generative AI and AI-assisted technologies in the manuscript preparation process}

During the revision of this manuscript, ChatGPT (OpenAI) was used to assist with checking mathematical arguments, clarifying definitions and proofs, organizing the exposition, and language editing. The authors take responsibility for the content of the manuscript.